\documentclass[final]{iopart}

\newcommand{\lnq}[1]{\ln_{q}(#1)}
\newcommand{\expq}[1]{\emph{e}_{q}(#1)}

\usepackage{iopams}  
\usepackage{amsopn}
\usepackage{array}  
\usepackage{tikz,pgfplots,amssymb}
\usepackage{etexcmds}
\usepackage{stackrel}
\usetikzlibrary{calc}
\usepgfplotslibrary{fillbetween}
\usetikzlibrary{arrows,intersections}

\DeclareMathOperator{\qplus}{\oplus_\mathit{q}}
\DeclareMathOperator{\qminus}{\ominus_\mathit{q}}
\DeclareMathOperator{\qtimes}{\otimes_\mathit{q}}
\DeclareMathOperator{\qdiv}{\oslash_\mathit{q}}
\DeclareMathOperator{\qproduct}{\odot_\mathit{q}}

\begin{document}
\title[q-Calculus Revisited]{\emph{q}-Calculus Revisited}

\author{Si Hyung Joo}

\address{Department of Industrial Engineering, Chonnam National University, 77, Yongbong-ro, Buk-gu, Gwangju, 61186, Republic of Korea}
\ead{innovation@jnu.ac.kr}
\vspace{10pt}
\begin{indented}
\item[]June 2021
\end{indented}

\begin{abstract}
In this study, a new representation is obtained for \emph{q}-calculus, as proposed by Borges [Phyica A 340 (2004) 95], and a new dual \emph{q}-integral is suggested. 

\vspace{2pc}
\end{abstract}

%
%
%
%

\section{Introduction}

With a generalization of the Boltzmann-Gibbs entropy \cite{Tsallis88}, the \emph{q}-logarithm and \emph{q}-exponential functions were first proposed by Tsallis  \cite{Tsallis94}.

\begin{eqnarray}
\lnq x &\equiv \frac{x^{1-q}-1}{1-\emph{q}}
\qquad \qquad \qquad & (x>0) \\
\expq x &\equiv [1+(1-q)x]_{+}^{1/(1-q)}
&(x,q \in \mathbb{R}) \\
& \mbox{where } [A]_+ \equiv \max \{ A,0 \} \nonumber
\end{eqnarray}

A \emph{q}-calculus associated with non-extensive statistical mechanics and thermodynamics was developed by Borges in 2004 \cite{Borges04}. He developed a primal \emph{q}-derivative operator, $D_{(q)}$, for which the \emph{q}-exponential function is an eigenfunction, as the ordinary exponential function is the eigenfunction of the ordinary derivative operator. 

\begin{equation}
D_{(q)} \Big( \expq x \Big) = \expq x \label{eq1}
\end{equation}

A primal \emph{q}-integral operator, $I_{(q)}$, which is the inverse operator of the primal \emph{q}-derivative operator, was also developed. The primal \emph{q}-integral of a \emph{q}-exponential function is a \emph{q}-exponential function. 

\begin{equation}
I_{(q)} \Big( \expq x \Big) = \expq x + c \label{eq2}
\end{equation}

In general, the following relationships hold for the primal \emph{q}-derivative and \emph{q}-integral.

\begin{eqnarray}
D_{(q)} \Big( {I_{(q)}}_a^x \Big( f(t) \Big) \Big) = f(x)  \label{eq3} \\
{I_{(q)}}_a^x \Big( D_{(q)} \Big( F(t) \Big) \Big) = F(x) +c \nonumber \\
\mbox{where } {I_{(q)}}_a^x \Big( f(t) \Big) = F(x) - F(a), F(t) ={I_{(q)}} \Big( f(t) \Big) - c \nonumber
\end{eqnarray}

For the \emph{q}-logarithm function, which is the inverse function of the \emph{q}-exponential function, he developed a dual \emph{q}-derivative operator, $D^{(q)}$. The dual \emph{q}-derivative of the \emph{q}-logarithm function is $1/x$, analogous to the fact that the ordinary derivative operator on the logarithm function gives $1/x$. 

\begin{equation}
D^{(q)} \Big( \lnq x \Big) = \frac{1}{x} \label{eq4}
\end{equation}

Finally, he suggested a dual \emph{q}-integral operator, $I^{(q)}$, which is the inverse operator of the dual \emph{q}-derivative. However, the dual \emph{q}-integral of $1/x$ is found to be $\ln (x) - \frac{1-q}{x} + c$, and the following relationship does not hold:

\begin{equation}
I^{(q)} \Big( \frac{1}{x} \Big) = \lnq x + c \label{eq5}
\end{equation}

In general, the following relationship does not hold for the dual \emph{q}-derivative and \emph{q}-integral, and it is a significant weakness of the dual \emph{q}-calculus suggested by Borges.

\begin{eqnarray}
D^{(q)} \Big( {I^{(q)}}_a^x \Big( f(t) \Big) \Big) = f(x) \label{eq6}  \\
{I^{(q)}}_a^x \Big( D^{(q)} \Big( F(t) \Big) \Big) = F(x) +c \nonumber \\
\mbox{where } {I^{(q)}}_a^x \Big( f(t) \Big) = F(x) - F(a), F(t)={I^{(q)}} \Big( f(t) \Big) - c  \nonumber
\end{eqnarray}

To address this issue, a new representation of \emph{q}-calculus with a new dual \emph{q}-integral operator is proposed herein, which satisfies equation (\ref{eq1}) $\sim$ (\ref{eq6}) with a modification of ordinary addition to \emph{q}-difference ($\qminus$) or \emph{q}-addition ($\qplus$) in equations (\ref{eq5}) and (\ref{eq6}): 

\begin{eqnarray}
D_{[q]} \Big( F(x) \Big) & \equiv \lim_{y \rightarrow x} \frac{F(x) - F(y)}{\ln[E_{q}(x)] - \ln[E_{q}(y)]} \\
\stackrel[{[q]}] {\vphantom{q}} {I} ^{x}_{x_0} \Big( f(t) \Big) 
& \equiv \int^{x}_{t=x_0} f(t) \ \mathrm{d} u(t) \\
& \quad \mbox{where } u(t) = \ln \Big[ \vert 1+(1-q)t \vert ^\frac{1}{1-q} \Big] \nonumber \\
D^{[q]}F(x) & \equiv \lim_{y \rightarrow x} \frac{ \ln( \expq {F(x)}) - \ln( \expq {F(y)})}{x - y} \\
\stackrel[\vphantom{q}] {{[q]}} {I} ^{x}_{x_0} \Big( f(t) \Big) 
& \equiv \ln_{q} \left[ \exp \left( \int^{x}_{x_0} f(t) \mathrm{d} t \right) \right] 
\end{eqnarray}

The new representation of \emph{q}-calculus is based on the concept of primal and dual \emph{q}-tangent lines, analogous to ordinary tangent lines. The primal and dual \emph{q}-derivatives are defined as the slope of the primal and dual \emph{q}-tangent lines at each point on the curve $y=f(x)$; this is analogous to the fact that the ordinary derivative of a function is the slope of the tangent line at each point on the curve $y=f(x)$. The primal and dual \emph{q}-integrals are defined as the signed primal and dual \emph{q}-area between the curve $y=f(x)$ and the horizontal axis, as the ordinary integral of a function is the signed area between the curve $y=f(x)$ and the horizontal axis.

The remainder of this paper is organized as follows. Section 2 provides some background on \emph{q}-algebra, \emph{q}-calculus of Borges, and ordinary calculus. Section 3 involves the derivation of a new representation of \emph{q}-calculus and the new dual \emph{q}-integral operator. Section 4 presents the relationship between the primal and dual \emph{q}-derivatives and integrals.

\section{Review of \emph{q}-Algebra, \emph{q}-Calculus and Ordinary Calculus}

\subsection{q-Algebra}

\emph{q}-algebra was first proposed by Borges \cite{Borges04}.

\begin{eqnarray} 
x \qplus y &\equiv x + y + (1-q)xy \\
x \qminus y &\equiv \frac{x - y}{1 + (1-q)y} \\
x \qtimes y &\equiv [x^{1-q}+y^{1-q}-1]_{+}^{1/(1-q)}
\qquad (x, y >0) \\
x \qdiv y &\equiv [x^{1-q}-y^{1-q}+1]_{+}^{1/(1-q)}
\qquad (x, y >0) \\
x^{\otimes_{q}^{n}} &\equiv \underbrace{ x \qtimes x \qtimes x \qtimes \cdots \qtimes x}_{n \ times} \
= [nx^{1-q}-(n-1)]_{+}^{1/(1-q)}\\
n \qproduct x &\equiv \underbrace{ x \qplus x \qplus x \qplus \cdots \qplus x}_{n \ times}  \
= \frac{1}{1-q}\{[1+(1-q)x]^{n}-1\}
\end{eqnarray}

The properties of the \emph{q}-logarithm and the \emph{q}-exponential can be expressed as follows.

\begin{eqnarray}
\lnq {xy} &= \lnq x \qplus \lnq y, \qquad
	\expq x \expq y &= \expq {x \qplus y} \\
\lnq {x \qtimes y} &= \lnq x + \lnq y, \qquad
	\expq x \qtimes \expq y &= \expq {x + y}\\
\lnq {x/y} &= \lnq x \qminus \lnq y, \qquad 
	\expq x / \expq y &= \expq {x \qminus y} \\
\lnq {x \qdiv y} &= \lnq x - \lnq y, \qquad 
	\expq x \qdiv \expq y &= \expq {x - y}
\end{eqnarray}

\subsection{q-Calculus}

Borges \cite{Borges04} defined primal and dual \emph{q}-derivatives and \emph{q}-integrals as follows. \\

\begin{eqnarray}
D_{(q)} \Big( f(x) \Big) &\equiv \lim_{y \rightarrow x} \frac{f(x)-f(y)}{x \qminus y}=[1+(1-q)x] \frac{\mathrm{d}f(x)}{\mathrm{d}x} \label{eqp1} \\
I_{(q)} \Big( f(x) \Big) &\equiv \int \frac{f(x)}{1+(1-q)x}{\mathrm{d}x} \label{eqp3}\\
D^{(q)} \Big( f(x) \Big) &\equiv \lim_{y \rightarrow x} \frac{f(x) \qminus f(y)}{x - y} \
=\frac{1}{1+(1-q)f(x)} \frac{\mathrm{d}f(x)}{\mathrm{d}x} \label{eqp2}\\
I^{(q)} \Big( f(x) \Big) &\equiv \int [1+(1-q)f(x)]f(x){\mathrm{d}x}
\end{eqnarray}

\subsection{Ordinary Calculus}

\subsubsection{Derivative}

Let $f(x)$ be the ordinary derivative function of $F(x)$.

\begin{equation}
\frac{\mathrm{d}}{\mathrm{d}x}F(x) \equiv 
\lim_{t \rightarrow x} \frac{F(x)-F(t)}{x - t} =f(x) \\
\end{equation}

A derivative operation is a function that takes a function $F(x)$ as an argument and produces another function $f(x)$ as an output.

In ordinary calculus, $f(x_0)$, the value of $f(x)$ evaluated at $x_0$, represents the slope of the tangent line $\mathcal{T}(C, P): y=T(x;F(x), x_0)=k_{(C,P)}x+c_{(C,P)}$ at the point $P=(x, y)=(x_0, F(x_0))$ on the curve $C: y=F(x)$.\\


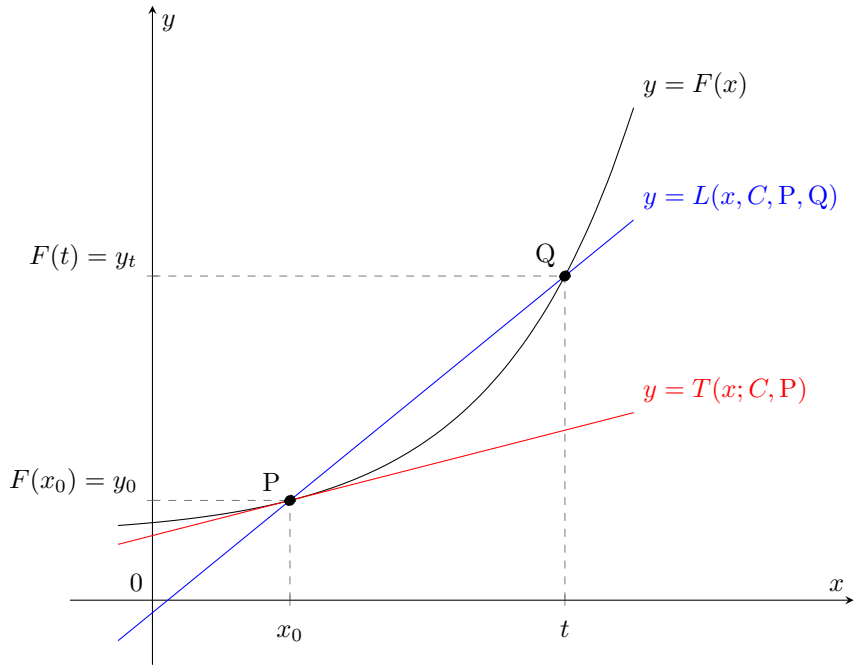
\begin{figure}[h]
\centering
\begin{tikzpicture}[
	declare function={f(\x)=0.2*exp(\x/2)+1;}
]

\newcommand\XTickA{2}
\newcommand\XTickB{6}

\begin{axis}[
  axis lines=center,
  xmin=-1.2,xmax=10.2,
  ymin=-1,ymax=9.2,
  width=12cm,
  xtick=data,
  ytick=data,
  typeset ticklabels with strut, 
  yticklabel style = {yshift=0.25cm},
  xticklabels={$x_0$,$t$},
  yticklabels={$F(x_0)=y_0$,$F(t)=y_t$},
  xlabel=$x$,
  ylabel=$y$,
  domain=-0.5:7,
  smooth
]

\addplot [ycomb,
          mark=*,
          mark options={black,mark size=2pt},
          gray, dashed,
          samples at={\XTickA,\XTickB}] {f(x)} 
          ;
\addplot [mark=none] coordinates {(\XTickA, {f(\XTickA)})} node[above left] {P};
\addplot [mark=none] coordinates {(\XTickB, {f(\XTickB)})} node[above left] {Q};

\addplot [mark=none] coordinates {(0,0)} node[above left] {0};

\addplot [xcomb, gray, dashed, samples at={\XTickA,\XTickB}] {f(x)};          

\addplot [name path = fcurve] {f(x)} node[above right] {$y=F(x)$};

\addplot [red, domain=-0.5:7] {(0.1*exp(\XTickA/2))*(x-\XTickA)+f(\XTickA)} node[above right] {$y=T(x;C,\mathrm{P})$};

\addplot [blue, domain=-0.5:7] {(f(\XTickB)-f(\XTickA))/(\XTickB-\XTickA)*(x-\XTickA)+f(\XTickA)} node[above right] {$y=L(x,C,\mathrm{P},\mathrm{Q})$};

\end{axis}
\end{tikzpicture}
\caption{Derivative as the slope of the tangent line} \label{fig:F1}
\end{figure}

Let the line passing through two points $P_i=(x_i,F(x_i))$ and $P_j=(x_j,F(x_j))$ on the curve $C$ be $\mathcal{L}(C, P_i, P_j): y=L(x; F(x), x_i, x_j)$.\\

\begin{equation}
L(x; F(x), x_i, x_j)= \frac{F(x_i)-F(x_j)}{x_i-x_j}(x-x_i)+F(x_i)
\end{equation}

The line $\mathcal{L}(C, P, Q)$ passing through $P=(x_0, F(x_0))$ and another point $Q=(t, F(t))$ on the curve $C$ tends to $\mathcal{T}(C, P)$, as $Q$ approaches $P$, and the slope of $\mathcal{L}(C, P, Q)$ converges to that of $ \mathcal{T}(C, P)$.

\begin{equation}
f(x_0) = \lim_{t \rightarrow x_0} \frac{L(x_0;F(x), x_0, t)-L(t;F(x), x_0, t)}{x_0 - t}=k_{(C,P)}
\end{equation}


Note that 
\begin{equation}
\frac{\mathrm{d}}{\mathrm{d}x} H(x) = \frac{\mathrm{d}}{\mathrm{d}x} F(x) \quad \Longleftrightarrow \quad H(x)=F(x)+c
\end{equation}

If $C$ is a line, that is, $C: F(x) = kx+c$, the curve (line) $C$ itself is regarded as the tangent line at each point on it, and the slope of the tangent line is constant, $k$, for all points on $C$.

\subsubsection{Integral}

A definite integral operation is a function that takes a function $f(x)$ and a pair of values $(x_L, x_H)$ as arguments and outputs the difference between the values of primitive function, $F(x)$, evaluated at $x_H$ and $x_L$.

\begin{equation}
\int_{x_L}^{x_H} f(x) \mathrm{d}x \equiv F(x_H) -F(x_L).
\end{equation}

Let the signed area between the curve $C':y=f(x)$ and the horizontal axis ($y=0$) from $x_L$ to $x_H$ be $A(f(x), x_L, x_H)$. The value of $A(f(x), x_L, x_H)$ is equal to the difference between the two values of the primitive function $F(x)$ evaluated at $x_H$ and $x_L$, \\

\begin{equation}
\int_{x_L}^{x_H} f(x) \mathrm{d}x \equiv F(x_H) -F(x_L) = A(f(x), x_L, x_H). 
\end{equation}

The relationship is clear when the primitive function is $F(x)=kx+c$, and its derivative is a constant; that is, $\frac{\mathrm{d}}{\mathrm{d}x}F(x)=k$ (see Figure 2).


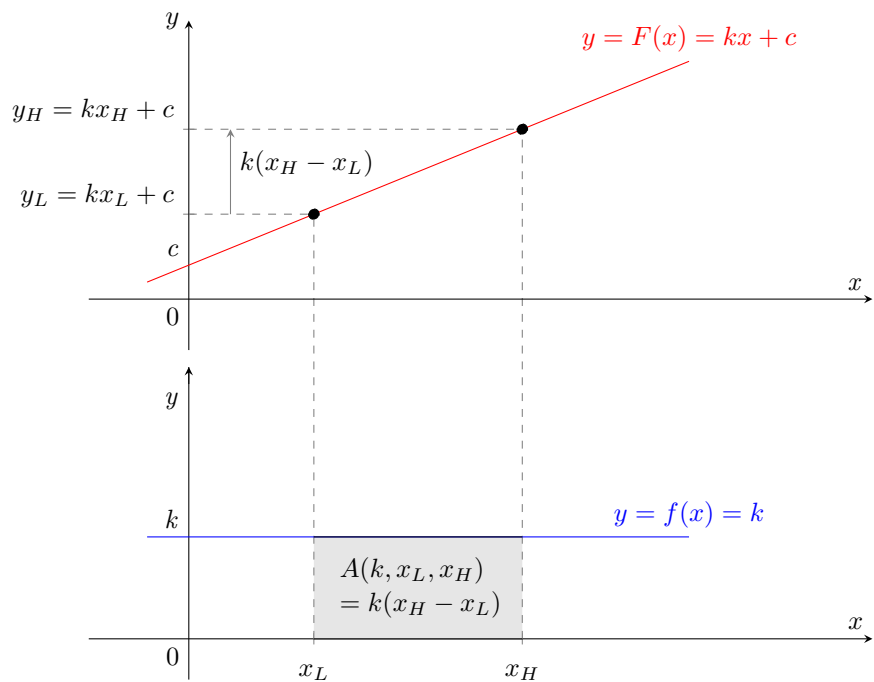
\begin{figure}[h]
\centering
\begin{tikzpicture}[ 	
	declare function={f(\x)=\x+11;},
	declare function={g(\x)= 3;},
	declare function={w(\x)= 10;}
]

\newcommand\XTickA{1.5}
\newcommand\XTickB{4}

\begin{axis}[
  axis lines=center,
  xmin=-1.2,xmax=8.2,
  ymin=-1.2,ymax=18.2,
  width=12cm,
  xtick=data,
  ytick=data,
  typeset ticklabels with strut, 
  yticklabel style = {yshift=0.25cm},
  xticklabels={$x_L$,$x_H$},
  yticklabels={$y_L=kx_L+c$,$y_H=kx_H+c$},
  xlabel=$x$,
  y label style={anchor=east},
  ylabel=$y$,
  domain=-0.5:6,
  smooth
]

\addplot [ycomb,
          mark=*,
          mark options={black,mark size=2pt},
          gray, dashed,
          samples at={\XTickA,\XTickB}] {f(x)};

\addplot [xcomb, gray, dashed, samples at={\XTickA,\XTickB}] {f(x)};          

\addplot [name path = fcurve, red] {f(x)} node[above] {$y=F(x)=kx+c$};

\addplot [mark=none] coordinates {(0,11)} node[above left] {\emph{c}};

\addplot [mark=none] coordinates {(0,10)} node[below left] {0};

\addplot [name path = gcurve, blue] {g(x)} node[above] {$y=f(x)=k$};

\addplot [mark=none] coordinates {(0,3)} node[above left] {\emph{k}};

\addplot [mark=none] coordinates {(0,0)} node[below left] {0};

\addplot [name path = g,domain=\XTickA:\XTickB] {g(x)};

\addplot[name path = h,domain=\XTickA:\XTickB] {0};

\addplot [mark=none] coordinates {(0,7)} node[left] {\emph{y}};

\addplot [mark=none, white, very thick] coordinates {(0, 8) (0, 8.5)};

\addplot [mark=none, -stealth, black] coordinates {(0, 7.5) (0, 8)};

\addplot [gray!20] fill between[of=g and h];

\addplot [name path = gcurve, -stealth, domain=-1.2:8.2] {w(x)} node[above left] {\emph{x}};

\addplot [mark=none, -stealth, gray] coordinates {(0.5, 12.5) (0.5, 15)};

\addplot [mark=none] coordinates {(0.5,14)} node[right] {$k(x_H-x_L)$};

\addplot [mark=none] coordinates {(1.7,2)} node[right] {$A(k,x_L,x_H)$};
\addplot [mark=none] coordinates {(1.7,1)} node[right] {$=k(x_H-x_L)$};

\end{axis}

\end{tikzpicture}

\caption{Definite integral and the signed area when $\frac{\mathrm{d}}{\mathrm{d}x}F(x)=k$} \label{fig:F2}
\end{figure}

\vspace{20pt}
In general, $A(f(x), x_L, x_H)$ can be evaluated as follows:


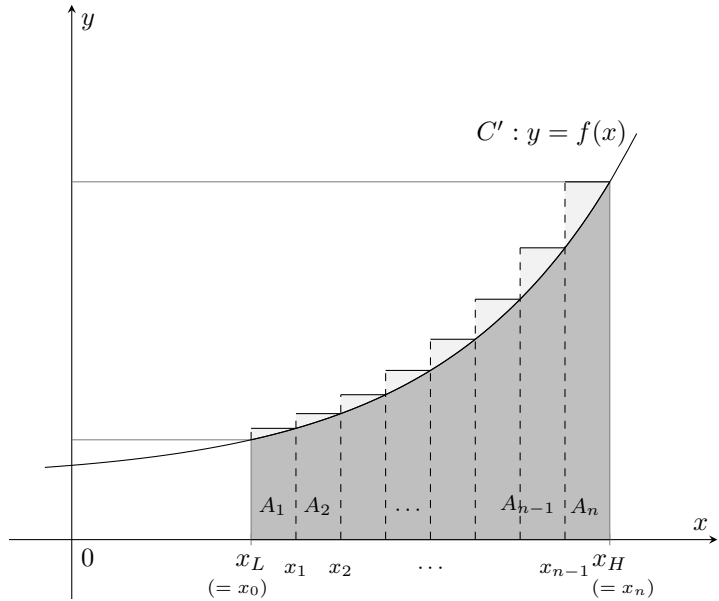
\begin{figure}[h]
\centering

\begin{tikzpicture}[ 	
	declare function={f(\x)=0.05*exp(\x/2)+0.2;}
]

\newcommand\XTickA{2}
\newcommand\XTickB{6}
\newcommand\NoP{8}

\begin{axis}[
  axis lines=center,
  xmin=-0.7,xmax=7.2,
  ymin=-0.2,ymax=1.8,
  width=11cm,
  xtick=data,
  ytick=\empty,
  xticklabels={$x_L$,$x_H$},
  xlabel=$x$,
  ylabel=$y$,
  domain=-0.3:6.3,
  smooth
]

\addplot [ycomb,
          mark=none,
          mark options={black,mark size=2pt},
          gray,
          samples at={\XTickA,\XTickB}] {f(x)};

\

\addplot [xcomb, gray, samples at={\XTickA,\XTickB}] {f(x)};          

\addplot [mark=none] coordinates {(0,0)} node[below right] {0};

\addplot [name path = fcurve] {f(x)} node[left] {$C':y=f(x)$};

\addplot [name path = f,domain=\XTickA:\XTickB] {f(x)};

\addplot[name path = g,domain=\XTickA:\XTickB] {0};

\foreach \t in {1,...,\NoP} {
	\addplot [name path = seg_f, domain=\XTickA+(\t-1)*((\XTickB-\XTickA)/\NoP):\XTickA+\t*((\XTickB-\XTickA)/\NoP)] {f(\XTickA+\t*((\XTickB-\XTickA)/\NoP))};
	\addplot [name path = seg_g, draw=none, domain=\XTickA+(\t-1)*((\XTickB-\XTickA)/\NoP):\XTickA+\t*((\XTickB-\XTickA)/\NoP)] {0};
	\addplot [gray!10] fill between [of = seg_f and seg_g];
}

	\addplot[dashed] coordinates{(\XTickA, {f(\XTickA)}) (\XTickA, {f(\XTickA+((\XTickB-\XTickA)/\NoP))})};

\addplot [gray!50] fill between[of=f and g];

\foreach \t in {2,...,\NoP} {
	\addplot[dashed] coordinates{(\XTickA+(\t-1)*((\XTickB-\XTickA)/\NoP), 0) (\XTickA+(\t-1)*((\XTickB-\XTickA)/\NoP), {f(\XTickA+\t*((\XTickB-\XTickA)/\NoP))})};
}

\addplot [mark=none] coordinates {({\XTickA-0.15},-0.1)} node[below] [font=\fontsize{7}{0}] {$(=x_0)$};
\addplot [mark=none] coordinates {({\XTickA+(1)*((\XTickB-\XTickA)/\NoP)},-0.05)} node[below] [font=\fontsize{8}{0}] {$x_1$};
\addplot [mark=none] coordinates {({\XTickA+(2)*((\XTickB-\XTickA)/\NoP)},-0.05)} node[below] [font=\fontsize{8}{0}] {$x_2$};
\addplot [mark=none] coordinates {({\XTickA+((floor(\NoP/2)+1)-1)*((\XTickB-\XTickA)/\NoP)},-0.05)} node[below] [font=\fontsize{8}{0}] {$\cdots$};
\addplot [mark=none] coordinates {({\XTickA+(\NoP-1)*((\XTickB-\XTickA)/\NoP)},-0.05)} node[below] [font=\fontsize{8}{0}] {$x_{n-1}$};
\addplot [mark=none] coordinates {({\XTickB+0.15},-0.1)} node[below] [font=\fontsize{7}{0}] {$(=x_n)$};

\addplot [mark=none] coordinates {({\XTickA+(1)*((\XTickB-\XTickA)/\NoP)},0.05)} node[above left] [font=\fontsize{8}{0}] {$A_1$};
\addplot [mark=none] coordinates {({\XTickA+(2)*((\XTickB-\XTickA)/\NoP)},0.05)} node[above left] [font=\fontsize{8}{0}] {$A_2$};
\addplot [mark=none] coordinates {({\XTickA+((floor(\NoP/2)+1)-1)*((\XTickB-\XTickA)/\NoP)},0.05)} node[above left] [font=\fontsize{8}{0}] {$\cdots$};
\addplot [mark=none] coordinates {({\XTickA+(\NoP-1)*((\XTickB-\XTickA)/\NoP)},0.05)} node[above left] [font=\fontsize{8}{0}] {$A_{n-1}$};
\addplot [mark=none] coordinates {({\XTickB},0.05)} node[above left] [font=\fontsize{8}{0}] {$A_n$};

\end{axis}

\end{tikzpicture}

\caption{General evaluation of signed area} \label{fig:F3}
\end{figure}

Let $T=\{x_0, x_1, \cdots, x_{n-1}, x_n\}$ be the $n$-partitions of $[x_L, x_H]$, $\Delta t_i = x_i-x_{i-1}$, and the length of the longest sub-interval is the norm of the partition, $\|T\|$. The signed area $A_i=A((f(x), x_{i-1}, x_{i})$ between the curve $C':y=f(x)$ and the horizontal axis ($y=0$) over the $i$-th partition $[x_{i-1}, x_i]$ can be approximated by the signed area of a rectangle with height $f(x_i)$ and base $\Delta t_i$; that is, $A_i \approx f(x_i)\Delta t_i$. Note that $A_i = f(x_i)\Delta t_i$ if $f(x)$ is constant over the partition.

The signed area between the curve $C'$ and the x axis over $[x_L, x_H]$ can be approximated by the sum of areas of rectangles, that is, $A(f(x), x_L, x_H) \approx \sum_{i=1}^{n} f(x_i)\Delta t_i$. Note that $A = \sum_{i=1}^{n} f(x_i)\Delta t_i$ if $f(x)$ is constant over $[x_L, x_H]$.    

As the norm of the partition approaches zero, the sum of areas of rectangles converges to the signed area between the curve $C$ and the horizontal axis ($y=0$) over $[x_L, x_H]$, and the definite integral is defined as the limit.

\begin{equation}
\int_{x_L}^{x_H} f(x) \mathrm{d}x = A(f(x), x_L, x_H)\equiv \lim_{\|T\| \rightarrow 0} \sum_{n=1}^{n} f(x_i)\Delta t_i
\end{equation}

In contrast, an indefinite integral is a function that takes a function $f(x)$ as an argument and produces a family of functions $\mathcal{F}= \left\{ F(x)+c, c \in \mathbb{R} \right\}$. The codomain of an indefinite integral operation is a set of families of functions, and the image of an indefinite integral is a translation family of functions along the y axis.

\begin{equation}
\int f(x)dx = F(x) + c = \mathcal{F}
\end{equation}

We can recover the exact function $F(x)$ only if we know a point $P=(x_0, F(x_0))$ on the curve $C$ with $f(x)$.

\begin{equation}
F(x) = \int_{x_0}^x f(t)dt + F(x_0) = A(f(x), x, x_0)+F(x_0)
\end{equation}

\newpage

\section{New Representation of \emph{q}-Calculus and New Dual \emph{q}-integral}

\subsection{Primal}

\subsubsection{Primal q-derivative}

Let a family of curves $\mathcal{L}_{q (k_q, \cdot)} = \left\{ y=L_{q}(x; k_q, c) \right\}$ have a constant primal \emph{q}-derivative (defined by equation (\ref{eqp1})), $k_q$,  at every point in their domain.

$L_{q}(x; k_q, c)$ must satisfy the condition $D_{(q)} \Big( L_{q}(x; k_q, c) \Big) = k_q$. 

\begin{equation}
D_{(q)} \Big( L_{q}(x; k_q, c) \Big) 
= \left\{ 1+(1-q)x \right\} \frac{\mathrm{d}}{\mathrm{d}x} L_{q}(x; k_q, c)
= k_q
\end{equation}

Therefore,
\begin{equation}
\eqalign{
L_{q}(x; k_q, c) & = \int \frac{k_q}{1+(1-q)x} \mathrm{d}x \cr
& = \frac{k_q}{1-q} \cdot \ln(\vert 1+(1-q)x \vert ) + c \cr
& = k_q \cdot \ln( \vert 1+(1-q)x \vert ^\frac{1}{1-q}]+c \cr
& = k_q \cdot \ln[ E_{q}(x) ] +c \\
} 
\end{equation}

where
\begin{equation}
E_{q}(x) \equiv \vert 1+(1-q)x \vert ^\frac{1}{1-q}
\end{equation}
Note that $E_{q}(x)= \expq x$, where $1+(1-q)x >0$.

Let us call the curve $\mathcal{L}_q (k_q, c):y=L_{q}(x; k_q, c)=k_q \cdot ln(E_{q}(x))+c$ a primal $\emph{q}$-line with primal $\emph{q}$-slope $k_q$ and $y$-intercept $c$.

Note that $\mathcal{L}_q (0, c):y=L_{q}(x; 0, c)$ is a horizontal line, $y=c$.

The line $L_{q}(x; k_q, c)$ has an interesting property:

\begin{equation}
\fl \frac{L_{q}(x_i; k_q, c)-L_{q}(x_j; k_q, c)}{\ln(E_{q}(x_i))-\ln(E_{q}(x_j))}=k_q, \forall x_i, x_j \in \mathbb{R}, x_i \neq x_j, x_i, x_j \neq \frac{1}{q-1} \label{eqpa1}
\end{equation}

From this finding, we define a new primal \emph{q}-derivative operator $D_{[q]}$, 

\begin{equation}
D_{[q]} \Big( F(x) \Big) \equiv \lim_{y \rightarrow x} \frac{F(x) - F(y)}{\ln[E_{q}(x)] - \ln[E_{q}(y)]}
\end{equation}

Note that 
\begin{equation}
D_{[q]} \Big( H(x) \Big) = D_{[q]} \Big( F(x) \Big) \quad \Longleftrightarrow \quad H(x)=F(x)+c.
\end{equation}

$\expq x$ is the eigenfunction of $D_{[q]}$.
\begin{equation}
\fl D_{[q]} \Big( \expq x \Big)
= \lim_{y \rightarrow x} \frac{\expq x - \expq y}{\ln[E_{q}(x)] - \ln[E_{q}(y)]} 
= \expq x, \ \mbox{if} \ \ 1+(1-q)x >0
\end{equation}

In general, $D_{[q]}$ is equivalent to $D_{(q)}$; that is, $D_{[q]} \Big( F(x) \Big) = D_{(q)} \Big( F(x) \Big)$, in general, because
\begin{equation}
\lim_{y \rightarrow x} \frac{\ln[E_{q}(x)] - \ln[E_{q}(y)]}{x \qminus y} = 1
\end{equation}

\begin{equation}
\eqalign{
D_{[q]} \Big( F(x) \Big) & \equiv \lim_{y \rightarrow x} \frac{F(x) - F(y)}{\ln[E_{q}(x)] - \ln[E_{q}(y)]} \cr
& = \lim_{y \rightarrow x} \frac{F(x) - F(y)}{\ln[E_{q}(x)] - \ln[E_{q}(y)]} \frac{\ln[E_{q}(x)] - \ln[E_{q}(y)]}{x \qminus y} \cr
& = \lim_{y \rightarrow x} \frac{F(x) - F(y)}{x \qminus y} = D_{(q)} \Big( F(x) \Big) \cr
}
\end{equation}

Let $\mathcal{L}_q(C, P_i, P_j) : y = L_q(x; F(x), x_i, x_j)$ be a primal $\emph{q}$-line that passes through the two points $P_i=(x_i, y_i)=(x_i, F(x_i))$ and $P_j=(x_j, y_j)=(x_i, F(x_i))$ on the curve $C: y=F(x)$.

$L_q(x; F(x), x_i, x_j))$ turns out to be
\begin{eqnarray}
L_q(x; F(x), x_i, x_j)) = k_q \cdot \ln(E_{q}(x))+c \\
k_q = \frac{F(x_i)-F(x_j)}{\ln(E_{q}(x_i))-\ln(E_{q}(x_j))} \\
c = y_i-k \cdot \ln(E_{q}(x_i))
\end{eqnarray}

$L_q(x; F(x), x_i, x_j))$ can also be expressed as follows:
\begin{eqnarray}
L_q(x; F(x), x_i, x_j)) = k_q \cdot \ln(E_{q}(x \qminus x_i))+ F(x_i) \\
L_q(x; F(x), x_i, x_j)) = k_q \cdot \ln(E_{q}(x \qminus x_j))+ F(x_j)
\end{eqnarray}


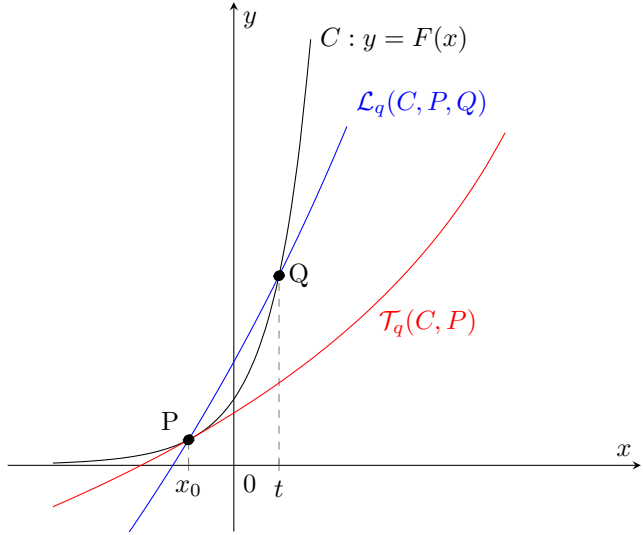
\begin{figure}[h]
\centering

\begin{tikzpicture}[ 	
	declare function={f(\x)= (max(1 - 0.1*x, 0))^(-10);}, 
	declare function={g(\x)= (1 + 0.1)^(-10)*(ln((max(1 - 0.1*(x+1), 0))^(-10))+1);}, 
	declare function={h(\x)= 
((1 + (-0.1))^(-10) - (1 + 0.1)^(-10))/(ln((1 + (-0.1))^(-10)) - ln((1 + 0.1)^(-10)))
*ln((max(1 - 0.1*x, 0))^(-10))
+(1 + 0.1)^(-10)
-((1 + (-0.1))^(-10) - (1 + 0.1)^(-10))/(ln((1 + (-0.1))^(-10)) - ln((1 + 0.1)^(-10)))
*ln((1 + 0.1)^(-10))
	;} 
]

\newcommand\XTickA{-1}
\newcommand\XTickB{1}

\begin{axis}[
  axis lines=center,
  xmin=-5,xmax=9,
  ymin=-1,ymax=7,
  width=10cm,
  xtick=data,
  ytick=\empty,
  yticklabel style = {above right},
  xticklabels={$x_0$,$t$},
  yticklabels={$F(x_0)$,$F(t)$},
  xlabel=$x$,
  ylabel=$y$,
  smooth
]

\addplot [ycomb,
          mark=*,
          mark options={black,mark size=2pt},
          gray, dashed,
          samples at={\XTickA,\XTickB}] {f(x)};


\addplot [mark=none] coordinates {(0,0)} node[below right] {0};

\addplot [name path = fcurve, domain=-4:1.7] {f(x)} node[right] {$C:y=F(x)$};
\addplot [red, name path = gcurve, domain=-4:6] {g(x)} node[below right]{};
\addplot [blue, name path = hcurve, domain=-4:2.5] {h(x)} node[above right] {$\mathcal{L}_q(C, P, Q)$};

\addplot [mark=none] coordinates {(\XTickA, {f(\XTickA)})} node[above left] {P};
\addplot [mark=none] coordinates {(\XTickB, {f(\XTickB)})} node[right] {Q};

\addplot [mark=none] coordinates {(3,2.5)} node[red, below right] {$\mathcal{T}_q(C, P)$};

\end{axis}
\end{tikzpicture}

\caption{Primal \emph{q}-derivative as the slope of the primal \emph{q}-tangent line} \label{fig:F4}
\end{figure}

We call $\mathcal{T}_q(C, P): y=T_q(x; F(x), x_0)=k_{q(C,P)} \cdot \ln(E_q(x))+c_{q(C,P)}$ as the primal $\emph{q}$-tangent line of the curve $C: y=F(x)$ at the point $P=(x_0, F(x_0))$ when $F(x_0)=T_{q}(x_0)$ and $D_{[q]} \Big( F(x) \Big) \Big\vert_{x=x_0} = k_{q(C,P)}$.

$D_{[q]} \Big( F(x) \Big) \Big\vert_{x=x_0}$ is the primal $\emph{q}$-slope of the primal $\emph{q}$-tangent line at point $P$ on the curve $C$.

The primal $\emph{q}$-line $\mathcal{L}_q(C, P, Q)$ passing through $P=(x_0, F(x_0))$ and another point $Q=(t, F(t))$ on the curve $C$ tends to $\mathcal{T}_q(C, P)$ as $Q$ approaches $P$, and the primal $\emph{q}$-slope of $\mathcal{L}_q(C, P, Q)$ converges to that of $\mathcal{T}_q(C, P)$.

\begin{eqnarray}
D_{[q]} \Big( F(x) \Big) \Big\vert_{x=x_0} & = \lim_{t \rightarrow x_0} \frac{L_q(x_0; F(x), x_0, t)-L_q(t; F(x), x_0, t)}{x_0 - t} \\
& =k_{q(C,P)} \nonumber
\end{eqnarray}

\subsubsection{Primal q-integral} If $f(x)= D_{[q]} \Big( F(x) \Big)$, the definite primal \emph{q}-integral of $f(x)$ from $x_L$ to $x_h$, denoted as $\stackrel[{[q]}] {\vphantom{q}} {I} ^{x_H}_{x_L} \Big( f(x) \Big)$,  is equal to $F(x_H) - F(x_L)$.

\begin{equation}
\stackrel[{[q]}] {\vphantom{q}} {I} ^{x_H}_{x_L} \Big( f(x) \Big) = F(x_H) - F(x_L)
\end{equation}

Equation (\ref{eqpa1}) gives us a hint on how a primal \emph{q}-integral can be related to a (deformed) signed area.

Consider the most immediate example, the case of primal \emph{q}-lines.
Let $F(x)=k_q \cdot ln(E_{q}(x))+c$ and $f(x)= D_{[q]} \Big( F(x) \Big) = k_q$. Let us consider the transformation $(x, y)=(x, f(x)) \rightarrow (u, v)= (\ln(E_{q}(x)), f(x))$.


\begin{figure}[h]
\centering

\begin{tikzpicture}[ 	
	declare function={f(\x)=(-1.5)*ln(abs(1-0.2*x))+11;}, 
	declare function={g(\x)= 3;},
	declare function={w(\x)= 10;}
]

\newcommand\XTickA{1}
\newcommand\XTickB{3.5}
\newcommand\XTickC{6.5}
\newcommand\XTickD{9}

\begin{axis}[
  axis lines=center,
  xmin=-1.2,xmax=10.2,
  ymin=-1.2,ymax=18.2,
  width=12cm,
  xtick=data,
  ytick=data,
  typeset ticklabels with strut, 
  xticklabel style = {yshift=0.25cm, scale=0.65},
  yticklabel style = {yshift=0.3cm, scale=0.85},
  xticklabels={$u_{L_1}$,$u_{H_1}$, $u_{H_2}$,$u_{L_2}$},
  yticklabels={${y_{L_{1}}, y_{H_{2}}}$,${y_{H_{1}}, y_{L_{2}}}$, ,},
  x label style={anchor=north, scale=0.8},
  xlabel=$u$,
  y label style={anchor=east, scale=0.85},
  ylabel=$y$,
  domain=-0.5:9,
  smooth
]

\addplot [ycomb,
          mark=*,
          mark options={black,mark size=1pt},
          gray, dashed,
          samples at={\XTickA,\XTickB,\XTickC,\XTickD}] {f(x)};
          
\addplot [mark=none, white, very thick] coordinates {(\XTickA, 7) (\XTickA, 10)};
\addplot [mark=none, white, very thick] coordinates {(\XTickB, 7) (\XTickB, 10)};
\addplot [mark=none, white, very thick] coordinates {(\XTickC, 7) (\XTickC, 10)};
\addplot [mark=none, white, very thick] coordinates {(\XTickD, 7) (\XTickD, 10)};

\addplot [xcomb, gray, dashed, samples at={\XTickC,\XTickD}] {f(x)};          

\addplot [name path = fcurve, red, domain=-0.5:4.8] {f(x)} node[above, scale=0.65] {$y=L_{q}(x; k_q,c)=k_{q} \cdot ln(E_{q}(x))+c$};

\addplot [name path = fcurve, red, domain=5.2:10] {f(x)} node[above] {};

\addplot [name path = hcurve, -stealth, domain=-1.2:10.2] {w(x)} node[below left, scale=0.8] {\emph{x}};
\addplot [mark=none] coordinates {(0,10)} node[below left, scale=0.85] {0};
\addplot [mark=none] coordinates {(0,11)} node[below left, scale=0.85] {\emph{c}};

\addplot [mark=none, purple, dashed] coordinates {(5, 9.5) (5, 16)};
\addplot [mark=none] coordinates {(5,9)} node[purple, scale=0.65] {$x=-\frac{1}{1-q}$};

\addplot [mark=none] coordinates {(\XTickA,9.3)} node[scale=0.8] {$x_{L_1}$};
\addplot [mark=none] coordinates {(\XTickB,9.3)} node[scale=0.8] {$x_{H_1}$};
\addplot [mark=none] coordinates {(\XTickC,9.3)} node[scale=0.8] {$x_{L_2}$};
\addplot [mark=none] coordinates {(\XTickD,9.3)} node[scale=0.8] {$x_{H_2}$};

\addplot [mark=none, dashed, gray] coordinates {(6.5, {(-1.5)*ln(abs(1-0.2*6.5))+11}) (7.3, {(-1.5)*ln(abs(1-0.2*6.5))+11})};

\addplot [mark=none, -stealth, gray] coordinates {(7, {(-1.5)*ln(abs(1-0.2*9))+11}) (7, {(-1.5)*ln(abs(1-0.2*6.5))+11})};

\addplot [mark=none] coordinates {(7.3,{(-1.5)*ln(abs(1-0.2*3.8))+11})} node[right, scale=0.75]  {$y_{H_{1}}-y_{L_{1}}$};

\addplot [mark=none] coordinates {(7.7,{(-1.5)*ln(abs(1-0.2*3))+11})} node[right, scale=0.75] {$= - (y_{H_{2}}-y_{L_{2}}) $};

\addplot [mark=none, white, very thick] coordinates {(0, 8) (0, 8.8)};
\addplot [mark=none, -stealth, black] coordinates {(0, 7.5) (0, 8)};
\addplot [mark=none] coordinates {(0,7)} node[left, scale=0.85] {\emph{v}};
\addplot [name path = gleft, blue,  domain= 0:4.5] {g(x)} node[above, scale=0.65] {};
\addplot [mark=none] coordinates {(0,3)} node[above left, scale=0.65] {$k_{q}$};
\addplot [mark=none] coordinates {(0,0)} node[below left, scale=0.65] {0};
\addplot [mark=none] coordinates {(4.7,0)} node[below left, scale=0.85] {\emph{u}};
\addplot [mark=none, white, very thick] coordinates {(4.2, 0) (4.8, 0)};
\addplot [mark=none, -stealth, black] coordinates {(4,0) (4.5,0)};

\addplot [mark=none] coordinates {(\XTickA,-0.8)} node[scale=0.5] {$=ln(E_q(x_{L_1}))$};
\addplot [mark=none] coordinates {(\XTickB,-0.8)} node[scale=0.5] {$=ln(E_q(x_{H_1}))$};

\addplot [mark=none] coordinates {(5,7)} node[left, scale=0.85] {\emph{v}};
\addplot [name path = gright, blue,  domain=5:9.5] {g(x)} node[above, scale=0.65] {};
\addplot [mark=none] coordinates {(5,3)} node[above left, scale=0.65] {$k_{q}$};
\addplot [mark=none] coordinates {(5,0)} node[below left, scale=0.65] {0};
\addplot [mark=none, -stealth, black] coordinates {(5, -1.2) (5, 8)};

\addplot [mark=none] coordinates {(\XTickC,-0.8)} node[scale=0.5] {$=ln(E_q(x_{H_2}))$};
\addplot [mark=none] coordinates {(\XTickD,-0.8)} node[scale=0.5] {$=ln(E_q(x_{L_2}))$};

\addplot [name path = g1,domain=\XTickA:\XTickB] {g(x)};
\addplot [name path = h1,domain=\XTickA:\XTickB] {0};

\addplot [name path = g2,domain=\XTickC:\XTickD] {g(x)};
\addplot [name path = h2,domain=\XTickC:\XTickD] {0};

\addplot [gray!10] fill between[of=g1 and h1];
\addplot [gray!40] fill between[of=g2 and h2];

\addplot [mark=none] coordinates {(1.2,1.5)} node[right, scale=0.65]  {$A_q(k_q, u_{L_1}, u_{H_1})$};

\addplot [mark=none] coordinates {(6.7,1.5)} node[right, scale=0.65]  {$A_q(k_q, u_{L_2}, u_{H_2})$};

\end{axis}

\end{tikzpicture}

\caption{Relationship between the primal \emph{q}-integral and signed primal \emph{q}-area} \label{fig:F5}
\end{figure}
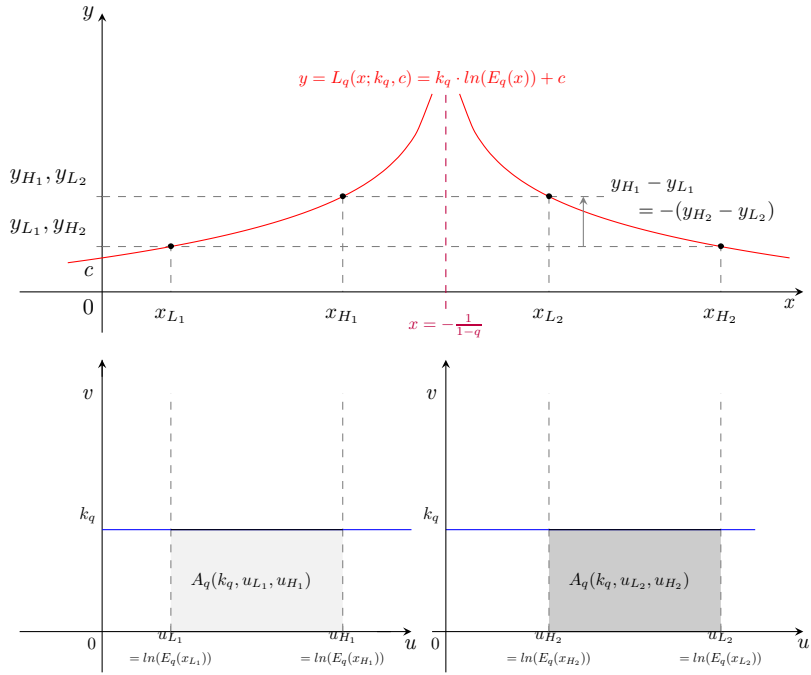
Figure 5 shows the graph of $y=F(x)$ at the top and that of $u$ and $v$ at the bottom. 
The bottom left graph shows the case $x_{L_1},x_{H_1}<-\frac{1}{1-q}$, and the bottom right graph shows the case $x_{L_2},x_{H_2}>-\frac{1}{1-q}$. $x_{L_2}$ and $x_{H_2}$ are set as $x_{L_2} = -\frac{2}{1-q} - x_{H_1}$ and $x_{H_2} = -\frac{2}{1-q} - x_{L_1}$.

${y_{L_{1}} = y_{H_{2}}}$ and ${y_{H_{1}} = y_{L_{2}}}$ because 
\begin{eqnarray}
E_{q}(x) & = \vert 1+(1-q)x \vert ^\frac{1}{1-q} \\ 
& = \vert -1-(1-q)x \vert ^\frac{1}{1-q} = \left \vert 1+(1-q) \left( - \frac{2}{1-q}-x  \right) \right \vert ^\frac{1}{1-q} \nonumber \\
& = E_{q}( - \frac{2}{1-q}-x ) \nonumber
\end{eqnarray}

Let $A_q(k_q, u_{L}, u_{H})$ be the signed rectangular area formed by $v=k_q$ and $v=0$ over the range $[u_L, u_H]$; this is called the primal \emph{q}-area.

\begin{equation}
A_q(k_q, u_L, u_H) = k_q \cdot ( u_H - u_L )
\end{equation}

With equation (\ref{eqpa1}), we find that the signed area of each shaded rectangle in the $u-v$ graph is equal to the corresponding definite integral.
\begin{eqnarray}
A_q(k_q, u_{L_1}, u_{H_1}) & = k_q \cdot ( u_{H_1} - u_{L_1} ) \\
& = k_q \cdot \left\{ \ln(E_{q}(x_{H_1}))-\ln(E_{q}(x_{L_1})) \right\} \nonumber \\
& = L_{q}(x_{H_1}; k_q, c) - L_{q}(x_{L_1}; k_q, c) \nonumber \\
& = \stackrel[{[q]}] {\vphantom{q}} {I} ^{x_{H_1}}_{x_{L_1}} \Big( L_{q}(x; k_q, c) \Big) \nonumber \\
A_q(k_q, u_{L_2}, u_{H_2}) & = k_q \cdot ( u_{H_2} - u_{L_2} ) \\
& = k_q \cdot \left\{ \ln(E_{q}(x_{H_2}))-\ln(E_{q}(x_{L_2})) \right\} \nonumber \\
& = L_{q}(x_{H_2}; k_q, c) - L_{q}(x_{L_2}; k_q, c) \nonumber \\
& = \stackrel[{[q]}] {\vphantom{q}} {I} ^{x_{H_2}}_{x_{L_2}} \Big( L_{q}(x; k_q, c) \Big) \nonumber
\end{eqnarray}

$y_{H_1} - y_{L_1} = -( y_{H_2} - y_{L_2} ), $ because ${y_{L_{1}} = y_{H_{2}}}$ and ${y_{H_{1}} = y_{L_{2}}}$. 
From this, we find that
\begin{equation}
\stackrel[{[q]}] {\vphantom{q}} {I} ^{x_H}_{x_L} \Big( L_{q}(x; k_q, c) \Big) 
= - \stackrel[{[q]}] {\vphantom{q}} {I} ^{-\frac{2}{1-q} - x_L}_{-\frac{2}{1-q} - x_H} \Big( L_{q}(x; k_q, c) \Big) 
\end{equation}

From this, we get  
\begin{eqnarray}
\fl \stackrel[{[q]}] {\vphantom{q}} {I} ^{x_H}_{x_L} \Big( L_{q}(x; k_q, c) \Big) 
& = \stackrel[{[q]}] {\vphantom{q}} {I} ^{-\frac{1}{1-q}}_{x_L} \Big( L_{q}(x; k_q, c) \Big) 
+\stackrel[{[q]}] {\vphantom{q}} {I} ^{x_H}_{-\frac{1}{1-q}} \Big( L_{q}(x; k_q, c) \Big) \\
& = \stackrel[{[q]}] {\vphantom{q}} {I} ^{-\frac{1}{1-q}}_{x_L} \Big( L_{q}(x; k_q, c) \Big)
- \stackrel[{[q]}] {\vphantom{q}} {I} ^{-\frac{1}{1-q}}_{-\frac{2}{1-q} - x_H} \Big( L_{q}(x; k_q, c) \Big) \nonumber \\ 
& = \stackrel[{[q]}] {\vphantom{q}} {I} ^{-\frac{2}{1-q} - x_H}_{x_L} \Big( L_{q}(x; k_q, c) \Big), \nonumber \\ 
& \quad \mbox{when} \ x_L < -\frac{1}{1-q} < x_H \nonumber
\end{eqnarray}

Based on the above findings, we define the definite primal \emph{q}-integral of $f(x)$ from $x_L$ to $x_h$, $\stackrel[{[q]}] {\vphantom{q}} {I} ^{x_H}_{x_L} \Big( f(x) \Big)$,  as follows:

\begin{eqnarray}
\stackrel[{[q]}] {\vphantom{q}} {I} ^{x_H}_{x_L} \Big( f(x) \Big) = F(x_H) - F(x_L) \equiv A_q(f(x), u_L, u_H)
\end{eqnarray}


We derive the primal \emph{q}-integral of the \emph{q}-exponential function, $\stackrel[{[q]}] {\vphantom{q}} {I} ^{x_H}_{x_L} \Big( \expq x \Big)$, with $A_q(k_q, u_{L}, u_{H})$. 

We can use $A_q(k_q, u_{L}, u_{H}) = k_q \cdot \left\{ \ln(\expq {x_{H}})-\ln(\expq {x_{L}}) \right\}$ instead of $A_q(k_q, u_{L}, u_{H}) = k_q \cdot \left\{ \ln(E_{q}(x_{H}))-\ln(E_{q}(x_{L})) \right\}$ because $k_q = \expq x =0$ when $1+(1-q)x \leq 0$.

Split $[x_L, x_H]$ into $n$-partitions, $T_q=\{x_L=x_0, x_1, x_2, \cdots, x_{n-1}, x_n = x_H\}$, where $x_i = x_L \qplus (i \qproduct t)$, 

\begin{equation}
t= \frac{1}{1-q} \left[ \left\{ 1 + (1-q) \cdot (x_H \qminus x_L) \right\}^{\frac{1}{n}}-1 \right]
\end{equation}


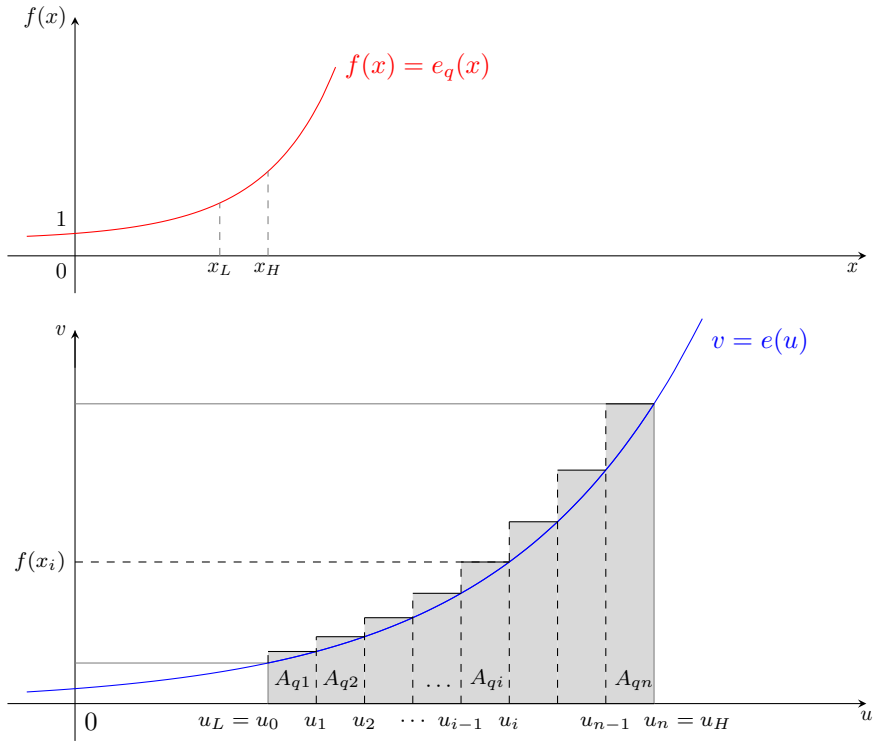
\begin{figure}[h]
\centering

\begin{tikzpicture}[ 	
	declare function={f(\x)=0.2*exp(\x/2);},
	declare function={g(\x)= 0.1*(max(1 - 0.1*x, 0))^(-10)+6.2;}, 
	declare function={w(\x)= 6;}
]

\newcommand\XTickA{2}
\newcommand\XTickB{6}
\newcommand\NoP{8}

\begin{axis}[
  axis lines=center,
  xmin=-0.7,xmax=8.2,
  ymin=-0.5,ymax=9.2,
  width=13cm,
  xtick=data,
  ytick=\empty,
  xtick=\empty,
  x label style={anchor=north, scale=0.8},
  xlabel=$u$,
  y label style={anchor=east, scale=0.85},
  ylabel=$f(x)$,
  domain=-0.5:7,
  smooth
]

\addplot [ycomb,
          mark=none,
          mark options={black,mark size=2pt},
          gray,
          samples at={\XTickA,\XTickB}] {f(x)};

\

\addplot [xcomb, gray, samples at={\XTickA,\XTickB}] {f(x)};          

\addplot [mark=none] coordinates {(0,0)} node[below right] {0};

\addplot [name path = ocurve, red, domain=-0.5:2.7] {g(x)} node[right] {$f(x)=e_q (x)$};
\addplot [name path = hcurve, -stealth, domain=-0.7:8.2] {w(x)} node[below left, scale=0.8] {\emph{x}};
\addplot [mark=none] coordinates {(0,6)} node[below left, scale=0.85] {0};
\addplot [mark=none] coordinates {(0,6.3)} node[above left, scale=0.85] {1};

\addplot [mark=none, white, very thick] coordinates {(0, 5) (0, 5.5)};
\addplot [mark=none, -stealth, black] coordinates {(0, 4.5) (0, 5)};
\addplot [mark=none] coordinates {(0,5)} node[left, scale=0.85] {\emph{v}};

\addplot [mark=none, dashed, gray] coordinates {(1.5, 6) (1.5, {0.1*(max(1 - 0.1*1.5, 0))^(-10)+6.2})};
\addplot [mark=none, dashed, gray] coordinates {(2, 6) (2, {0.1*(max(1 - 0.1*2, 0))^(-10)+6.2})};
\addplot [mark=none] coordinates {(1.5,6)} node[below, scale=0.8] {$x_{L}$};
\addplot [mark=none] coordinates {(2,6)} node[below, scale=0.8] {$x_{H}$};

\addplot [name path = fcurve, blue, domain=-0.5:6.5] {f(x)} node[below right] {$v=e(u)$};

\addplot [name path = f, blue, domain=\XTickA:\XTickB] {f(x)};

\addplot[name path = g,domain=\XTickA:\XTickB] {0};

\foreach \t in {1,...,\NoP} {
	\addplot [name path = seg_f, domain=\XTickA+(\t-1)*((\XTickB-\XTickA)/\NoP):\XTickA+\t*((\XTickB-\XTickA)/\NoP)] {f(\XTickA+\t*((\XTickB-\XTickA)/\NoP))};
	\addplot [name path = seg_g, draw=none, domain=\XTickA+(\t-1)*((\XTickB-\XTickA)/\NoP):\XTickA+\t*((\XTickB-\XTickA)/\NoP)] {0};
	\addplot [gray!30] fill between [of = seg_f and seg_g];
}

	\addplot[dashed] coordinates{(\XTickA, {f(\XTickA)}) (\XTickA, {f(\XTickA+((\XTickB-\XTickA)/\NoP))})};

\foreach \t in {2,...,\NoP} {
	\addplot[dashed] coordinates{(\XTickA+(\t-1)*((\XTickB-\XTickA)/\NoP), 0) (\XTickA+(\t-1)*((\XTickB-\XTickA)/\NoP), {f(\XTickA+\t*((\XTickB-\XTickA)/\NoP))})};
}

\addplot [mark=none] coordinates {({\XTickA-0.3},-0.05)} node[below] [font=\fontsize{8}{0}] {$u_L=u_0$};
\addplot [mark=none] coordinates {({\XTickA+(1)*((\XTickB-\XTickA)/\NoP)},-0.05)} node[below] [font=\fontsize{8}{0}] {$u_1$};
\addplot [mark=none] coordinates {({\XTickA+(2)*((\XTickB-\XTickA)/\NoP)},-0.05)} node[below] [font=\fontsize{8}{0}] {$u_2$};
\addplot [mark=none] coordinates {({\XTickA+((floor(\NoP/2))-1)*((\XTickB-\XTickA)/\NoP)},-0.05)} node[below] [font=\fontsize{8}{0}] {$\cdots$};
\addplot [mark=none] coordinates {({\XTickA+((floor(\NoP/2)+1)-1)*((\XTickB-\XTickA)/\NoP)},-0.05)} node[below] [font=\fontsize{8}{0}] {$u_{i-1}$};
\addplot [mark=none] coordinates {({\XTickA+(\NoP-3)*((\XTickB-\XTickA)/\NoP)},-0.05)} node[below][font=\fontsize{8}{0}] {$u_{i}$};
\addplot [mark=none] coordinates {({\XTickA+(\NoP-1)*((\XTickB-\XTickA)/\NoP)},-0.05)} node[below] [font=\fontsize{8}{0}] {$u_{n-1}$};
\addplot [mark=none] coordinates {({\XTickB+0.35},-0.05)} node[below] [font=\fontsize{8}{0}] {$u_n=u_H$};

\addplot [mark=none] coordinates {({\XTickA+(1)*((\XTickB-\XTickA)/\NoP)+0.05},0.05)} node[above left] [font=\fontsize{8}{0}] {$A_{q1}$};
\addplot [mark=none] coordinates {({\XTickA+(2)*((\XTickB-\XTickA)/\NoP)+0.05},0.05)} node[above left] [font=\fontsize{8}{0}] {$A_{q2}$};
\addplot [mark=none] coordinates {({\XTickA+((floor(\NoP/2)+1)-1)*((\XTickB-\XTickA)/\NoP)},0.05)} node[above left] [font=\fontsize{8}{0}] {$\cdots$};
\addplot [mark=none] coordinates {({\XTickA+(5)*((\XTickB-\XTickA)/\NoP)+0.05},0.05)} node[above left] [font=\fontsize{8}{0}] {$A_{qi}$};
\addplot [mark=none] coordinates {({\XTickB+0.1},0.05)} node[above left] [font=\fontsize{8}{0}] {$A_{qn}$};

\addplot[dashed] coordinates{(0, {f(\XTickA+5*((\XTickB-\XTickA)/\NoP))}) ({\XTickA+(5)*((\XTickB-\XTickA)/\NoP)}, {f(\XTickA+5*((\XTickB-\XTickA)/\NoP))})};

\addplot[mark=none] coordinates {(0, {f(\XTickA+5*((\XTickB-\XTickA)/\NoP))})} node[left] [font=\fontsize{8}{0}] {$f(x_i)$};

\end{axis}

\end{tikzpicture}

\caption{Definite primal \emph{q}-integral of a \emph{q}-exponential function} \label{fig:F6}
\end{figure}

$T_q$ transforms into $U_q = \left\{ \ln [ \expq {x_L} ] = u_0, u_1, u_2, \cdot, u_{n-1}, u_n = \ln [ \expq {x_H} ] \right\}$, where $u_i = \ln[ \expq {x _i} ]$, with a transformation $(x,y) = (x, f(x))=(x, \expq x) \rightarrow (u, v) = (\ln[ \expq x ],f(x)) = (\ln[ \expq x ], \expq x )$.

\begin{equation}
u_i = \ln [ \expq {x_i} ] = \ln [ \expq {x_L} ] + i \cdot \ln \left[ \left\{ \expq t \right\} \right]
\end{equation}

\begin{equation}
u_i - u_{i-1} = \ln \left[ \left\{ \expq t \right\} \right]
\end{equation}

\begin{eqnarray}
\expq t & = \left\{1 + (1-q) \cdot (x_H \qminus x_L) \right\}^{\frac{1}{n(1-q)}} \\
& = \left\{ \expq {x_H \qminus x_L} \right\}^{\frac{1}{n}} \nonumber \\
\ln \Big( \expq t \Big) & = \frac{1}{n(1-q)} \ln \left[ 1 + (1-q) \cdot (x_H \qminus x_L) \right] \label{eqpa2}\\
& =  \frac{1}{n} \ln \left[ \expq {x_H \qminus x_L} \right] \nonumber
\end{eqnarray}

Equation (\ref{eqpa2}) shows that the norm of the partition $U_q$, $\|U_q\|$, approaches zero when the number of partitions approaches infinity.

Let $A_{qi}$ be the primal \emph{q}-area formed by $v = \expq {x_i}$ and $v=0$ over the $i$-th partition $[u_{i-1}, u_i]$.

\begin{equation}
\fl \expq {x_i} = \expq {x_L \qplus (i \qproduct t)} = \expq {x_L} \cdot \expq {i \qproduct t} = \expq {x_L} \cdot \left\{ \expq t \right\} ^i
\end{equation}

\begin{eqnarray}
A_{qi} & = A_q(\expq {x_i}, u_{i-1}, u_i) = \expq {x_i} \cdot (u_i - u_{i-1}) \\
& = \left[ \expq {x_L} \cdot \left\{ \expq t \right\} ^i \right] \ln \left[ \expq t \right] \nonumber \\
& = \expq {x_L} \cdot \left\{ 1 + (1-q) \cdot (x_H \qminus x_L) \right\}^{\frac{i}{n(1-q)}} \cdot \frac{1}{n(1-q)} \ln \left[ 1 + (1-q) \cdot (x_H \qminus x_L) \right] \nonumber \\
& = \expq {x_L} \cdot \frac{\ln(z)}{1-q} \cdot \frac{1}{n} {\left( z^{\frac{1}{n(1-q)}} \right) }^i  \nonumber \\
& \quad \mbox{where} \ z = 1 + (1-q) \cdot (x_H \qminus x_L) \nonumber
\end{eqnarray}

$A_q(k_q, u_{L}, u_{H})$ can be evaluated as the sum of $A_{qi}$ when the norm of the partition $U_q$, $\|U_q\|$, approaches zero,
\begin{eqnarray}
\sum_{i=1}^{n}  A_{qi} & = \sum_{i=1}^{n} \left\{ \expq {x_L} \cdot \frac{\ln(z)}{1-q} \cdot \frac{1}{n} {\left( z^{\frac{1}{n(1-q)}} \right) }^i \right\} \\
& = \expq {x_L} \cdot \frac{\ln(z)}{1-q}  \cdot \frac{1}{n} \cdot \frac{z^{\frac{1}{n(1-q)}} \left(z^{\frac{n}{n(1-q)}}-1 \right) }{z^{\frac{1}{(1-q)}}-1} \nonumber \\
& = \expq {x_L} \cdot \frac{\ln(z)}{1-q}  \cdot \left(z^{\frac{1}{1-q}}-1 \right)  \cdot \frac{1}{n}  \cdot \frac{z^{\frac{1}{n(1-q)}}}{z^{\frac{1}{n(1-q)}}-1} \nonumber
\end{eqnarray}

\begin{eqnarray}
A_q(k_q, u_{L}, u_{H}) &= \lim_{\|U_q\| \rightarrow 0} \sum_{i=1}^{n}  A_{(q)i} \\
& = \expq {x_L} \cdot \frac{\ln(z)}{1-q} \cdot \left(z^{\frac{1}{1-q}}-1 \right) \cdot \lim_{n \rightarrow \infty } \Big( \frac{1}{n}  \frac{z^{\frac{1}{n(1-q)}}}{z^{\frac{1}{n(1-q)}}-1} \Big) \nonumber \\
& = \expq {x_L} \cdot \frac{\ln(z)}{1-q} \cdot \left(z^{\frac{1}{1-q}}-1 \right) \cdot \frac{1-q}{\ln (z)} \nonumber \\
& = \expq {x_L} \cdot \left( z^{\frac{1}{1-q}}-1 \right) \nonumber \\
& = \expq {x_L} \cdot \left\{ \expq {x_H \qminus x_L} - 1 \right\} \nonumber \\
& = \expq {x_L} \cdot \left\{ \frac{\expq {x_H}}{\expq {x_L}} -1 \right\} \nonumber \\
& = \expq {x_H} - \expq {x_L} \nonumber
\end{eqnarray}

Therefore,

\begin{equation}
\stackrel[{[q]}] {\vphantom{q}} {I} ^{x_H}_{x_L} \Big( \expq x \Big) = F(x_H) -F(x_L) = \expq {x_H} - \expq {x_L}
\end{equation}

If $F(x) = \expq x$, we know that $D_{[q]} \Big( F(x) \Big) = \expq x$ and $F(0) = 1$, for all \emph{q}. 
\begin{eqnarray}
F(x) = 1 + \stackrel[{[q]}] {\vphantom{q}} {I} ^{x}_{0} \Big( \expq t \Big) =  1 + \left\{ \expq {x} - \expq {0} \right\} = \expq {x}
\end{eqnarray}

Indefinite primal \emph{q}-integral of the \emph{q}-exponential function can be expressed as follows:

\begin{eqnarray}
F(x) = \stackrel[{[q]}] {\vphantom{q}} {I} \Big( \expq x \Big) = \stackrel[{[q]}] {\vphantom{q}} {I} ^{x}_{x_0} \Big( \expq t \Big) =  \expq {x} - \expq {x_0} = \expq {x} + c
\end{eqnarray}

In general, the signed area between the trajectory formed by $(u, v)$ over the range $[x_L, x_H]$ and $v=0$ can be evaluated using the Riemann–Stieltjes integral, $\int^{x_H}_{x=x_L} f(x) \ \mathrm{d} u(x)$ \cite{Lang93}, where $u(x) = \ln(E_{q}(x))$.

Therefore, when $f(x)= D_{[q]} \Big( F(x) \Big)$,
\begin{eqnarray}
\stackrel[{[q]}] {\vphantom{q}} {I} ^{x_H}_{x_L} \Big( f(x) \Big) & = F(x_H) - F(x_L) \label{eqpa3} \\
& = \int^{x_H}_{x=x_L} f(x) \ \mathrm{d} u(x) \nonumber \\
& = \int^{x_H}_{x=x_L} \frac{f(x)}{1+(1-q)x}  \ \mathrm{d} x \nonumber
\end{eqnarray}

Equation (\ref{eqpa3}) shows that the new primal \emph{q}-integral, $\stackrel[{[q]}] {\vphantom{q}} {I}$, is equivalent to the primal \emph{q}-integral $I_{(q)}$ in equation (\ref{eqp3}).

Equation (\ref{eq5}) holds for $D_{[q]}$ and $\stackrel[{[q]}] {\vphantom{q}} {I}$ because they are equivalent to $D_{(q)}$ and $I_{(q)}$, respectively.

When $f(x)= D_{[q]} \Big( F(x) \Big)$ and $F(x_0) = -c$, the new primal indefinite \emph{q}-integral can be expressed as follows:

\begin{equation}
\stackrel[{[q]}] {\vphantom{q}} {I} ^{x}_{x_0} \Big( f(t) \Big) = F(x) + c
\end{equation}


\newpage

\subsection{Dual}

\subsubsection{Dual q-derivative}

Let a family of curves $\mathcal{L}^{q}_{(k^q, \cdot)} = \left\{ y=L^{q}(x; k^q, c) \right\}$ have a constant dual \emph{q}-derivative (as defined by equation (\ref{eqp2})), $k^q$,  at every point on their domain.

$L^{q}(x; k^q, c)$ must satisfy the condition $D^{(q)} L^{q}(x; k^q, c) = k^q$, 

\begin{equation}
D^{(q)} L^{q}(x; k^q, c) 
= \frac{1}{1+(1-q)L^{q}(x; k^q, c)} \frac{\mathrm{d}}{\mathrm{d}x} L^{q}(x; k^q, c)
= k^q
\end{equation}

Therefore,
\begin{equation}
\eqalign{
L^{q} (x; k^q, c) & = c' \cdot \exp((1-q)k^q x) - \frac{1}{1-q} \cr
& = \frac{c}{1-q} \exp((1-q)k^q x)- \frac{1}{1-q} \cr
& = \left\{ \frac{1}{1-q} \exp((1-q)k^q x) + \frac{1}{1-q} \right\} + \frac{c-1}{1-q} \cr
& \quad + (1-q) \frac{c-1}{1-q} \left\{ \frac{1}{1-q} \exp((1-q)k^q x) + \frac{c-1}{1-q} \right\} \cr
& = \left\{ \frac{1}{1-q} \exp((1-q)k^q x) + \frac{1}{1-q} \right\} \qplus \frac{c-1}{1-q} \cr
& = \lnq {\exp(k^q x)} \qplus \frac{c-1}{1-q} 
} 
\end{equation}

Note that ${\mathcal{L}^q}_{(k^q, \cdot)}$ is a $\qplus$ translation family of a curve along the y axis. 

Consider the curve $\mathcal{L}^q(k^q, c):y=L^{q}(x; k^q, c)= \lnq { \exp(k^q x)} \qplus \frac{c-1}{1-q}$ as the dual $\emph{q}$-line with dual $\emph{q}$-slope $k^q$ and $y$-intercept $\frac{c-1}{1-q}$.

$L^{q}(x; k^q, c)$ also has the property,

\begin{equation}
\fl \frac{ \ln( \expq {L^{q}(x_i; k^q, c)}) - \ln( \expq {L^{q}(x_j; k^q, c)})}{x_i-x_j}=k^q, \forall x_i, x_j \in \mathbb{R}, x_i \neq x_j .\label{eqda4}
\end{equation}

From this finding, we define a new dual \emph{q}-derivative operator $D^{[q]}$ as 

\begin{equation}
D^{[q]} \Big( F(x) \Big) \equiv \lim_{y \rightarrow x} \frac{ \ln( \expq {F(x)}) - \ln( \expq {F(y)})}{x - y}
\end{equation}

Note that 
\begin{equation}
D^{[q]} \Big( H(x) \Big) = D^{[q]} \Big( F(x) \Big) \quad \Longleftrightarrow \quad H(x)=F(x) \qplus c \label{eqda1}
\end{equation}

It follows that $D^{[q]} \lnq x = \frac{1}{x}$.
\begin{equation}
\eqalign{
D^{[q]} \lnq x 
& = \lim_{y \rightarrow x} \frac{ \ln( \expq {\lnq x}) - \ln( \expq {\lnq y})}{x - y} \cr
& = \lim_{y \rightarrow x} \frac{ \ln(x) - \ln(y)}{x - y} = \frac{1}{x}
}
\end{equation}

$D^{[q]}$ turns out to be equivalent to $D^{(q)}$; that is, $D^{[q]}F(x) = D^{(q)}F(x)$, in general, because
\begin{equation}
\lim_{y \rightarrow x} \frac{ \ln(\expq {F(x)}) - \ln(\expq {F(y)})}{F(x) \qminus F(y)} = 1
\end{equation}

\begin{equation}
\eqalign{
D^{[q]} \Big( F(x) \Big) & \equiv 
\lim_{y \rightarrow x} \frac{ \ln( \expq {F(x)}) - \ln( \expq {F(y)})}{x - y} \cr
& = \lim_{y \rightarrow x} \frac{ \ln( \expq {F(x)}) - \ln( \expq {F(y)})}{F(x) \qminus F(y)} 
\frac{F(x) \qminus F(y)}{x - y} \cr
& = \lim_{y \rightarrow x} \frac{F(x) \qminus F(y)}{x - y} = D^{(q)}\Big( F(x) \Big) \cr
}
\end{equation}


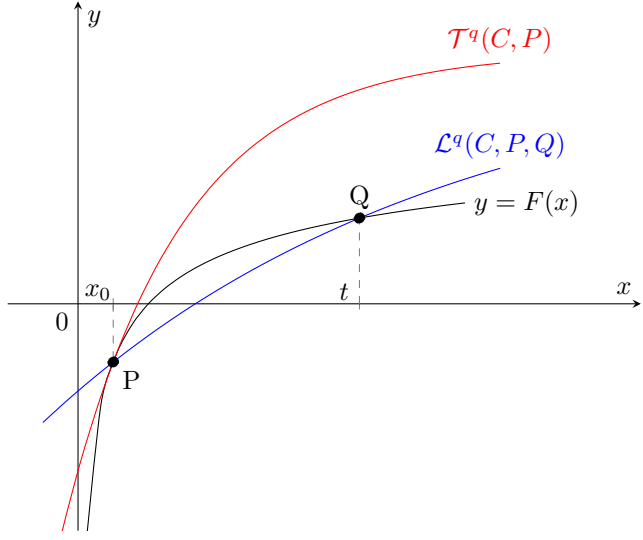
\begin{figure}[h]
\centering

\begin{tikzpicture}[ 	
	declare function={f(\x)= (x^(-0.3)-1)/(-0.3);}, 
	declare function={g(\x)= (((0.5)/exp(1)*exp(2*x))^(-0.3)-1)/(-0.3);}, 
	declare function={h(\x)=
	(((((0.5)/exp((0.5)*(1/(3.5)*ln((1+(-0.3)*(((4*2)^(-0.3)-1)/(-0.3)))^(1/(-0.3))))))^(-0.3))^(1/(-0.3))
*exp((1/(3.5)*ln((1+(-0.3)*(((4*2)^(-0.3)-1)/(-0.3)))^(1/(-0.3))))*x))^(-0.3)-1)/(-0.3)
	;} 
]

\newcommand\XTickA{0.5}
\newcommand\XTickB{4}

\begin{axis}[
  axis lines=center,
  xmin=-1,xmax=8,
  ymin=-3,ymax=4,
  width=10cm,
  xtick=data,
  ytick=\empty,
  xticklabel style = {xshift=-0.2cm, above},
  yticklabel style = {above right},
  xticklabels={$x_0$,$t$},
  yticklabels={$F(x_0)$,$F(t)$},
  xlabel=$x$,
  ylabel=$y$,
  smooth
]

\addplot [ycomb,
          mark=*,
          mark options={black,mark size=2pt},
          gray, dashed,
          samples at={\XTickA,\XTickB}] {f(x)};


\addplot [mark=none] coordinates {(0,0)} node[below left] {0};

\addplot [name path = fcurve, domain=0.1:5.5] {f(x)} node[right] {$y=F(x)$};
\addplot [red, name path = gcurve, domain=-0.5:6] {g(x)} node[above] {$\mathcal{T}^{q}(C, P)$};
\addplot [blue, name path = hcurve, domain=-0.5:6] {h(x)} node[above] {$\mathcal{L}^{q}(C, P, Q)$};

\addplot [mark=none] coordinates {(\XTickA, {f(\XTickA)})} node[below right] {P};
\addplot [mark=none] coordinates {(\XTickB, {f(\XTickB)})} node[above] {Q};


\end{axis}
\end{tikzpicture}

\caption{Dual \emph{q}-derivative as the slope of the dual \emph{q}-tangent line} \label{fig:F7}
\end{figure}

Let $\mathcal{L}^q(C, P_i, P_j): y = L^q(x; F(x), x_i, x_j)$ be the dual $\emph{q}$-line that passes through the two points $P_i=(x_i, y_i)=(x_i, F(x_i))$ and $P_j=(x_j, y_j)=(x_j, F(x_j))$ on the curve $C: y=F(x)$.

$L^{q}(x; F(x), P_i, P_j))$ turns out to be

\begin{equation}
L^{q}(x; F(x), x_i, x_j)) 
= \lnq { \exp(k^q x)} \qplus \frac{c-1}{1-q} 
\end{equation}
\begin{eqnarray}
k^q = \frac{ \ln( \expq {F(x_i)}) - \ln( \expq {F(x_j)})}{x_i-x_j} \\
\frac{c-1}{1-q} = y_i \qminus \lnq { \exp(k^q x_i)}
\end{eqnarray}

$L^{q}(x; F(x), x_i, x_j))$ can also be expressed as follows:
\begin{eqnarray}
L^{q}(x; F(x), x_i, x_j)) = \lnq {\exp(k^q (x-x_i)} \qplus F(x_i)  \\
L^{q}(x; F(x), x_i, x_j)) = \lnq {\exp(k^q (x-x_j)} \qplus F(x_j) 
\end{eqnarray}

We call $\mathcal{T}^{q}(C, P): y=T^{q}(x; F(x), x_0)=\lnq {\exp({k^q}_{(C, P)} x)} \qplus \frac{c^q_{(C,P)}-1}{1-q}$ as the dual $\emph{q}$-tangent line of the curve $C: y=F(x)$ at the point $P=(x_0, F(x_0))$ when $F(x_0)=T^{q}(x_0)$ and $D^{[q]} \Big( F(x) \Big) \Big\vert_{x=x_0} = {k^q}_{(C,P)}$.

$D^{[q]} \Big( F(x) \Big) \Big\vert_{x=x_0}$ is the dual $\emph{q}$-slope of the dual $\emph{q}$-tangent line at the point $P$ on the curve $C$.

The dual $\emph{q}$-line $L^{q}(x; C, P, Q)$ passing through $P=(x_0, F(x_0))$ and another point $Q=(t, F(t))$ on the curve $C$ tends to $\mathcal{T}^{q}(C, P)$ as $Q$ approaches $P$, and the dual $\emph{q}$-slope of $L^{q}(x; C, P, Q)$ converges to that of $\mathcal{T}^{q}(C, P)$.

\begin{eqnarray}
\fl D^{[q]} \Big( F(x) \Big) \Big\vert_{x=x_0} & = \lim_{t \rightarrow x_0} \frac{ \ln(\expq {L^{q}(x_0; F(x), x_0, t)}) - \ln(\expq {L^{q}(t; F(x), x_0, t)})}{x_0 - t} \\ 
& = {k^q}_{(C,P)} \nonumber
\end{eqnarray}

\subsubsection{New dual q-integral}
Let the indefinite dual \emph{q}-integral operator, $\stackrel[\vphantom{q}] {{[q]}} {I}$, be the inverse operator of the dual \emph{q}-derivative operator, $D^{[q]}$.
Because the $\qplus$-translation family of a function has the same \emph{q}-derivative (see equation (\ref{eqda1})), the indefinite dual \emph{q}-integral operator should produce a $\qplus$-translation family of a function.

\begin{equation}
f(x)= D^{[q]} \Big( F(x) \Big) \quad \Longrightarrow \quad \stackrel[\vphantom{q}] {{[q]}} {I} \Big( f(x) \Big) = F(x) \qplus c.
\end{equation}

Let us denote the definite dual \emph{q}-integral of $f(x)$ from $x_L$ to $x_h$ as $\stackrel[\vphantom{q}] {{[q]}} {I} ^{x_H}_{x_L} \Big( f(x) \Big)$.

It is desirable and natural that a definite dual \emph{q}-integral vanishes when the integrating range is zero.

\begin{equation}
\stackrel[\vphantom{q}] {{[q]}} {I} ^{c}_{c} \Big( f(x) \Big) = 0,\qquad \forall c  \in \mathbb{R} . \label{eqda2}
\end{equation}

From equation (\ref{eqda2}), the function $F(x)$, with $D^{[q]} \Big( F(x) \Big) = f(x)$ and $F(x_0)=0$, can be represented as follows:

\begin{equation}
F(x) = \stackrel[\vphantom{q}] {{[q]}} {I} ^{x}_{x_0} \Big( f(t) \Big).
\end{equation}

If $H(x)= F(x) \qplus c$, $D^{[q]} \Big( H(x) \Big) = f(x)$ and $H(x_0)= F(x_0) \qplus c = c$. $H(x)$ can be represented as 
\begin{eqnarray}
H(x) & = F(x) \qplus c  = \stackrel[\vphantom{q}] {{[q]}} {I} ^{x}_{x_0} \Big( f(t) \Big) \qplus c \\
& = \stackrel[\vphantom{q}] {{[q]}} {I} ^{x}_{x_0} \Big( f(t) \Big) \qplus H(x_0) \nonumber
\end{eqnarray}

Therefore,
\begin{eqnarray}
\stackrel[\vphantom{q}] {{[q]}} {I} ^{x}_{x_0} \Big( f(t) \Big) & = H(x) \qminus H(x_0) = \Big( F(x) \qplus c \Big) \qminus \Big( F(x_0) \qplus c \Big) \\
& = F(x) \qminus F(x_0) \nonumber
\end{eqnarray}

Therefore, 

\begin{eqnarray}
\stackrel[\vphantom{q}] {{[q]}} {I} ^{x_H}_{x_L} \Big( f(x) \Big) & = F(x_H) \qminus F(x_L) = \lnq { \frac{\expq {F(x_H)}}{\expq {F(x_L)}} } \label{eqda3}
\end{eqnarray}

From equation (\ref{eqda3}), we find that the following relationship holds for the definite dual \emph{q}-integral.

\begin{equation}
\stackrel[\vphantom{q}] {{[q]}} {I} ^{x_H}_{x_L} \Big( f(x) \Big) 
= \stackrel[\vphantom{q}] {{[q]}} {I} ^{c}_{x_L} \Big( f(x) \Big) 
\qplus \stackrel[\vphantom{q}] {{[q]}} {I} ^{x_H}_{c} \Big( f(x) \Big),\qquad \forall c  \in \mathbb{R}
\end{equation}

Equation (\ref{eqda4}) gives us a hint on how a definite dual \emph{q}-integral can be related to a signed area.

Consider the immediate example, the case of dual \emph{q}-lines.
Let $F(x)=L^{q}(x; k^q, c)=\lnq {\exp(k^q x)} \qplus \frac{c-1}{1-q}$ and $f(x)= D^{[q]} \Big( F(x) \Big) = k^q$. 

\begin{eqnarray}
\stackrel[\vphantom{q}] {{[q]}} {I} ^{x_H}_{x_L} \Big( k^q \Big) 
& = L^{q}(x_H; k^q, c) \qminus L^{q}(x_L; k^q, c) \\
& = \lnq { \frac{\expq {L^{q}(x_H; k^q, c)}}{\expq {L^{q}(x_L; k^q, c)}} } \nonumber
\end{eqnarray}

Let us consider the transformation $(x,Y)=(x, F(x))=(x, L^{q}(x; k^q, c)) \rightarrow (x, w)= (x, \ln ( \expq Y )) = (x, \ln ( \expq {L^{q}(x; k^q, c)} )) $.


\begin{figure}[h]
\centering

\begin{tikzpicture}[ 	
	declare function={f(\x)=(-4)*exp(((-0.2)*2.5*\x))+5+21;}, 
	declare function={h(\x)=\x+11;},
	declare function={g(\x)= 3;}
]

\newcommand\XTickA{1}
\newcommand\XTickB{5}

\begin{axis}[
  axis lines=center,
  xmin=-1.2,xmax=10.2,
  ymin=-1.2,ymax=28.2,
  width=12cm,
  xtick=data,
  ytick=data,
  typeset ticklabels with strut, 
  yticklabel style = {yshift=0.3cm, scale=0.65},
  xticklabels={$x_L$,$x_H$},
  yticklabels={$Y_L$,$Y_H$},
  xlabel=$x$,
  y label style={anchor=east, scale=0.85},
  ylabel=$Y$,
  domain=-0.5:9,
  smooth
]

\addplot [ycomb,
          mark=*,
          mark options={black,mark size=1pt},
          gray, dashed,
          samples at={\XTickA,\XTickB}] {f(x)};

\addplot [xcomb, gray, dashed, samples at={\XTickA,\XTickB}] {f(x)};          

\addplot [name path = fcurve, red, domain=-0.15:7] {f(x)} node[above, scale=0.65] {$Y=L^{q}(x; k^q,c)= \lnq {exp(kx)} \qplus \frac{c-1}{1-q}$};

\addplot [name path = hcurve1, -stealth, domain=-1.2:10.2] {20} node[above left] {\emph{x}};
\addplot [mark=none] coordinates {(0,20)} node[below left, scale=0.65] {0};
\addplot [mark=none] coordinates {(0,22)} node[left] {$\frac{c-1}{1-q} $};

\addplot [name path = hcurve1, -stealth, domain=-1.2:10.2] {10} node[above left] {\emph{x}};
\addplot [mark=none] coordinates {(0,10)} node[below left, scale=0.65] {0};

\addplot [mark=none] coordinates {(0,17)} node[right, scale=0.65] {$w= \ln (\emph{e}_q(Y))$};
\addplot [mark=none] coordinates {(7.3,15)} node[right, blue, scale=0.65] {$=k^qx + \ln[ \expq {\frac{c-1}{1-q}} ]$};
\addplot [mark=none, white, very thick] coordinates {(0, 18) (0, 18.8)};
\addplot [mark=none, -stealth, black] coordinates {(0, 17.5) (0, 18)};

\addplot [name path = hcurve, blue, domain=-0.5:7] {h(x)} node[below right, scale=0.65] {$w= \ln \left[ \expq {L^{q}(x; k^q, c)} \right]$};

\addplot [mark=none, gray, dashed] coordinates {(0, 12) (1, 12)};
\addplot [mark=none, gray, dashed] coordinates {(0, 16) (5, 16)};

\addplot [mark=*, mark options={black,mark size=1pt}] coordinates {(1, 12)};
\addplot [mark=*, mark options={black,mark size=1pt}] coordinates {(5, 16)};

\addplot [mark=none] coordinates {(0, 12)} node[left, scale=0.65] {$w_L$};
\addplot [mark=none] coordinates {(0, 16)} node[left, scale=0.65] {$w_H$};

\addplot [mark=none, -stealth, gray] coordinates {(0.5, 12) (0.5, 16)};

\addplot [mark=none] coordinates {(0.5,14)} node[right, scale=0.75]  {$w_{H}-w_{L}$};

\addplot [mark=none] coordinates {(0,7)} node[left, scale=0.65] {\emph{y}};
\addplot [mark=none, white, very thick] coordinates {(0, 8) (0, 8.8)};
\addplot [mark=none, -stealth, black] coordinates {(0, 7.5) (0, 8)};

\addplot [mark=none] coordinates {(0,0)} node[below left, scale=0.65] {0};

\addplot [name path = gcurve, blue, domain=-0.5:7] {g(x)} node[above, scale=0.65] {$y=D^{[q]} L^{q}(x; k^q, c)$};

\addplot [mark=none] coordinates {(0,3)} node[above left, scale=0.65] {$k^{q}$};

\addplot [name path = g,domain=\XTickA:\XTickB] {g(x)};
\addplot [name path = h,domain=\XTickA:\XTickB] {0};
\addplot [gray!20] fill between[of=g and h];

\addplot [mark=none] coordinates {(1.5,1.5)} node[right, scale=1]  {$A^q(k^q, x_{L}, x_{H})$};

\end{axis}

\end{tikzpicture}

\caption{Relationship between the dual \emph{q}-integral and signed dual \emph{q}-area} \label{fig:F8}
\end{figure}
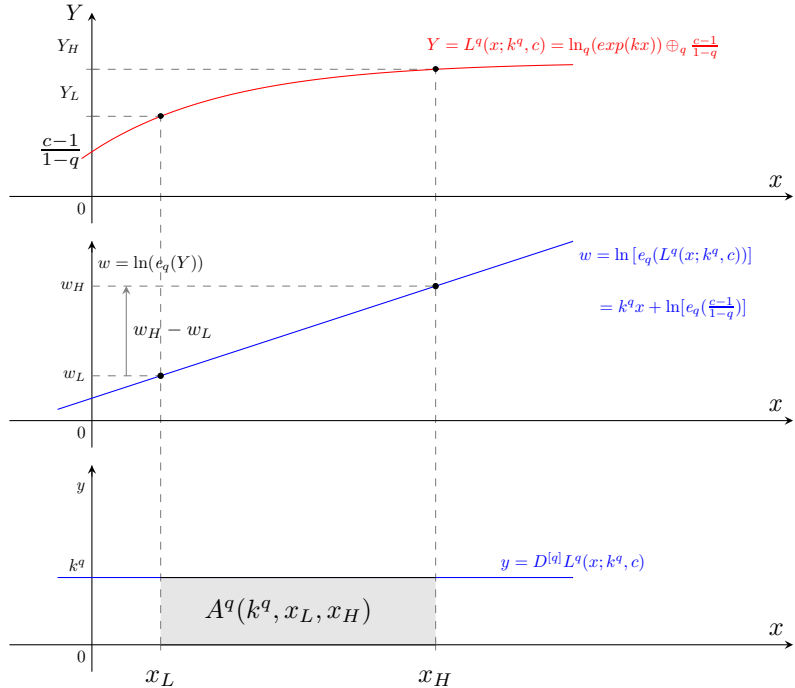

Figure 8 shows the graph of $Y=F(x)$ at the top, the graph of $x$ and $w$ in the middle, and the graph of $y=f(x)$ at the bottom. 

With equation (\ref{eqda4}), we find that the signed area of the shaded rectangle at the bottom graph is equal to $w_H-w_L$ in the middle graph.

Let $A(k^q, x_{L}, x_{H})$ be the ordinary signed rectangular area formed by $y=k^q$ and $y=0$ over the range $[x_L, x_H]$,

\begin{eqnarray}
A(k^q, x_L, x_H) & = k^q \cdot ( x_H - x_L ) \\
& = w_H-w_L \nonumber \\
& = \ln \left[ \expq {L^{q}(x_H; k^q, c)} \right] -\ln \left[ \expq {L^{q}(x_L; k^q, c)} \right] \nonumber \\
& = \ln \Big[ \frac{\expq {L^{q}(x_H; k^q, c)}}{\expq {L^{q}(x_L; k^q, c)}} \Big] \nonumber \\
& = \ln \Big[ \expq {\lnq {\frac{\expq {L^{q}(x_H; k^q, c)}}{\expq {L^{q}(x_L; k^q, c)}}}} \Big] \nonumber \\
& = \ln \left[ \expq {L^{q}(x_H; k^q, c) \qminus L^{q}(x_L; k^q, c)}  \right] \nonumber \\
& = \ln \left[ \expq {\stackrel[\vphantom{q}] {{[q]}} {I} ^{x_H}_{x_L} \Big( k^q \Big)}  \right] \nonumber
\end{eqnarray}

Let $A^q(k^q, x_L, x_H) \equiv \lnq {\exp \left[ A(k^q, x_L, x_H) \right]}=\lnq {\exp \left[ k^q \cdot ( x_H - x_L ) \right]}$ be a dual \emph{q}-area,

\begin{eqnarray}
\stackrel[\vphantom{q}] {{[q]}} {I} ^{x_H}_{x_L} \Big( k^q \Big) 
& = \lnq {\exp \left[ A(k^q, x_L, x_H) \right]} \\
& = A^q(k^q, x_L, x_H) \nonumber
\end{eqnarray}

For $f(x)=D^{[q]} \Big( F(x) \Big)$, in general, the ordinary signed area, $A(f(x), x_{L}, x_{H})$, formed by $y=f(x)$ and $y=0$ over the range $[x_L, x_H]$ is equal to the ordinary definite integral $\int^{x_H}_{x_L} f(x) \mathrm{d}x$. Therefore, the following relationship holds true:

\begin{eqnarray}
\stackrel[\vphantom{q}] {{[q]}} {I} ^{x_H}_{x_L} \Big( f(x) \Big) 
& = A^q(f(x), x_L, x_H) \label{eqda5}\\
& = \lnq {\exp \left[ A(f(x), x_L, x_H) \right]} \nonumber \\
& = \ln_q \Big( \exp \left[ \int^{x_H}_{x_L} f(x) \mathrm{d}x \right] \Big) \nonumber
\end{eqnarray}

Using equation (\ref{eqda5}), $\stackrel[\vphantom{q}] {{[q]}} {I} ^{x}_{x_0} \Big( \frac{1}{t} \Big)$ is found to be a \emph{q}-logarithmic function, 

\begin{eqnarray}
\stackrel[\vphantom{q}] {{[q]}} {I} ^{x}_{x_0} \Big( \frac{1}{t} \Big) 
& = \ln_q \Big( \exp \left[ \int^{x}_{x_0} \frac{1}{t} \mathrm{d}t \right] \Big) \\
& = \lnq {\exp \left[ \ln (x) - \ln (x_0) \right]} \nonumber \\
& = \lnq {\frac{x}{x_0} } = \lnq {x} \qminus \lnq {x_0} \nonumber
\end{eqnarray}

If $F(x) = \lnq x$, we know that $D^{[q]} \Big( \lnq x \Big) = \frac{1}{x}$ and $F(1) = 0$ for all \emph{q}, 
\begin{eqnarray}
F(x) = \stackrel[\vphantom{q}] {{[q]}} {I} ^{x}_{1} \Big( \frac{1}{t} \Big) =  \lnq {x} \qminus \lnq {1} = \lnq {x}
\end{eqnarray}

The indefinite dual \emph{q}-integral of $\frac{1}{x}$ can be expressed as follows, and equation (\ref{eq5}) holds with a modification of ordinary addition to \emph{q}-difference ($\qminus$).
\begin{eqnarray}
F(x) = \stackrel[\vphantom{q}] {{[q]}} {I} \Big( \frac{1}{t} \Big) 
\stackrel[\vphantom{q}] {{[q]}} {I} ^{x}_{x_0} \Big( \frac{1}{t} \Big)
= \lnq {x} \qminus \lnq {x_0} = \lnq {x} \qminus c.
\end{eqnarray}

Equation (\ref{eq6}) also holds with a modification of ordinary addition to \emph{q}-difference ($\qminus$).
\begin{eqnarray}
D^{[q]} \Big( \stackrel[\vphantom{q}] {{[q]}} {I} ^{x}_{a} \Big( f(x) \Big) \Big) & = D_{[q]} \Big( F(x) \qminus F(a) \Big) = f(x) \\
\stackrel[\vphantom{q}] {{[q]}} {I} ^{x}_{a} \Big( D^{[q]} \Big( F(x) \Big) \Big) & = \stackrel[\vphantom{q}] {{[q]}} {I} ^{x}_{a} \Big( f(x) \Big) = F(x) \qminus c
\end{eqnarray}

\section{Relationship Between the Primal and Dual}

\subsection{Relationship Between Primal and Dual q-derivatives}

Let $G(y)$ be the inverse function of $F(x)$, that is, $y=F(x)$ and $x=G(y)$, and $P=(x_0, y_0)$, $Q=(x_1, y_1)$ be points on the curve $C: y=F(x)$. The curve $C$ can also be represented by $C: x=G(y)$.

Let $\mathcal{L}_q(C, P, Q)$ be a primal \emph{q}-line that passes through $P$ and $Q$.

\begin{equation}
\mathcal{L}_q(C, P, Q): y = k_q \cdot \ln(E_{q}(x \qminus x_1))+ y_1 \label{eqpp1}
\end{equation}

$\mathcal{L}_q(C, P, Q)$ can also be represented as
\begin{equation}
\mathcal{L}_q(C, P, Q): x = \lnq {\exp(\frac{1}{k_q} (y-y_1)} \qplus x_1 \label{eqpp2}
\end{equation}

Equations (\ref{eqpp1}) and (\ref{eqpp2}) show that the primal \emph{q}-line that passes through $P$ and $Q$ is the dual \emph{q}-line passing through $P$ and $Q$, $\mathcal{L}^q(C, P, Q)$, and $k^q$; the dual \emph{q}-slope of $\mathcal{L}^q(C, P, Q)$ is $\frac{1}{k_q}$.

Let $\mathcal{T}_q(C, P)$ and $\mathcal{T}^q(C, P)$ be the primal and dual \emph{q}-tangent lines at $P$ on $C$, respectively, and let ${k_q}^*$ and $k^{q*}$ be the primal and dual \emph{q}-slopes of the corresponding \emph{q}-tangent line.
As $Q$ approaches $P$, $k_q$ converges to ${k_q}^*$, and $k^q$ converges to $k^{q*}$.
Therefore, the primal \emph{q}-derivative and the dual \emph{q}-derivative are inversely related.
\begin{equation}
k^{q*} = \frac{1}{{k_q}^*}.
\end{equation}

\subsection{Relationship Between Primal and Dual q-integrals}

Let $f(x)=D_{[q]} \Big( F(x) \Big)$ and $g(x)=D^{[q]} \Big( G(x) \Big)$,

\begin{eqnarray}
\stackrel[{[q]}] {\vphantom{q}} {I} ^{x_1}_{x_0} \Big( f(x) \Big) = y_1 - y_0 \qquad
\stackrel[\vphantom{q}] {{[q]}} {I} ^{y_1}_{y_0} \Big( g(y) \Big) = x_1 \qminus x_0 \\
\frac{\stackrel[{[q]}] {\vphantom{q}} {I} ^{x_1}_{x_0} \Big( f(x) \Big)}{\stackrel[\vphantom{q}] {{[q]}} {I} ^{y_1}_{y_0} \Big( g(y) \Big)} = \frac{y_1 - y_0}{x_1 \qminus x_0} = k_q = \frac{1}{k^q} \\
\qquad \mbox{where } k_q = \frac{y_1 - y_0}{\ln \left[ E_{q}(x_1)\right]-\ln \left[ E_{q}(x_0) \right]} \nonumber
\end{eqnarray}

Therefore, the definite primal and dual \emph{q}-integrals are proportionally related to the \emph{q}-slope of the \emph{q}-line passing through two points $P$ and $Q$.

\ack 
The author would like to thank Prof. Borges for his helpful advice and insights on the issue examined in this study.


\section*{Reference}
\begin{thebibliography}{5}
	\bibitem[1]{Tsallis88} Tsallis C 1988 Possible generalization of Boltzmann-Gibbs statistics \emph{Journal of Statistical Physics}, 52, 479-487
	\bibitem[2]{Tsallis94} Tsallis C 1994 What are the numbers that experiments provide? \emph{Quimica Nova}, 17, 468-471.
	\bibitem[3]{Borges04} Borges E P 2004 A possible deformed algebra and calculus inspired in nonextensive thermostatistics, \emph{Physica A}, 340, 95-101.
	\bibitem[4]{Lang93} Lang S 1993 \emph{Real and Functional Analysis} 3rd edn (New York: Springer)
\end{thebibliography}
\end{document}